\documentclass[fontsize=11pt]{scrartcl}

\usepackage{amsmath}
\usepackage{amssymb}
\usepackage{amsthm}
\usepackage{stmaryrd}
\usepackage{bbm}
\usepackage{mathtools}
\usepackage{mathrsfs}
\usepackage{color}
\usepackage{graphicx}
\usepackage{a4,fullpage}
\usepackage{enumerate}
\usepackage[normalem]{ulem}
\usepackage[hyperindex,breaklinks,hidelinks]{hyperref}

\usepackage{url}

\usepackage{color}

\usepackage{tikz}
\usetikzlibrary{decorations.pathmorphing,decorations.pathreplacing,calc,arrows.meta}

\newcommand{\bE}{\ensuremath{\mathbb{E}}}

\newcommand{\bN}{\ensuremath{\mathbb{N}}}

\newcommand{\bP}{\ensuremath{\mathbb{P}}}

\newcommand{\bR}{\ensuremath{\mathbb{R}}}

\newcommand{\bZ}{\ensuremath{\mathbb{Z}}}

\newcommand{\ind}[1]{\mathbbm{1}_{\{#1\}}} 
\newcommand{\indset}[1]{\mathbbm{1}_{#1}} 

\newcommand{\cC}{\ensuremath{\mathcal{C}}}
\newcommand{\cD}{\ensuremath{\mathcal{D}}}

\newcommand{\cF}{\ensuremath{\mathcal{F}}}

\newcommand{\cP}{\ensuremath{\mathcal{P}}}

\newcommand{\cR}{\ensuremath{\mathcal{R}}}

\newcommand{\cU}{\ensuremath{\mathcal{U}}}

\newcommand{\ann}{{\mathrm{ann}}}
\newcommand{\pre}{{\mathrm{pre}}}
\newcommand{\que}{{\mathrm{que}}}
\newcommand{\annpre}{{\mathrm{ann}\times\mathrm{pre}}}
\newcommand{\boxpre}{{\mathrm{box-que}\times\mathrm{pre}}}

\newcommand{\abs}[1]{\vert #1 \vert}
\newcommand{\Abs}[1]{\Bigl\vert #1 \Bigr\vert}

\newcommand{\norm}[1]{\left\Vert #1 \right\Vert}

\newcommand{\ddx}[1][1]{\ifnum#1=1 \frac{d}{dx} \else \frac{d^{#1}}{dx^{#1}} \fi}
\newcommand{\ddy}[1][1]{\ifnum#1=1 \frac{d}{dy} \else \frac{d^{#1}}{dy^{#1}} \fi}
\newcommand{\ddt}[1][1]{\ifnum#1=1 \frac{d}{dt} \else \frac{d^{#1}}{dt^{#1}} \fi}

\newcommand{\compl}{\mathsf{C}}

\usepackage{todonotes} 

\theoremstyle{plain}
\newtheorem{thm}{Theorem}[section]
\newtheorem{theorem}{Theorem}[section]  

\newtheorem{proposition}[theorem]{Proposition}

\newtheorem{lemma}[theorem]{Lemma}
\newtheorem{remark}[theorem]{Remark}

\theoremstyle{definition}

\newtheorem{definition}[thm]{Definition}

\theoremstyle{remark}

\numberwithin{equation}{section}

\newcommand{\simu}{\mathrm{sim}}

\newcommand{\wh}{\widehat}

\def\MR#1{\href{http://www.ams.org/mathscinet-getitem?mr=#1}{MR#1}}

\allowdisplaybreaks[4]

\title{Quenched local limit theorem for a directed random walk on the backbone of a supercritical oriented percolation cluster for $d \ge 1$}
\author{Stein Andreas Bethuelsen, Matthias Birkner, \\ Andrej
	Depperschmidt and Timo Schl\"{u}ter}
\date{\today}

\begin{document}
	\maketitle
	
	\abstract{In this work we extend the quenched local limit theorem of
		\cite{BethuelsenBirknerDepperschmidtSchluter}. More precisely, we
		consider a directed random walk on the backbone of the supercritical
		oriented percolation cluster in dimensions $d+1$ with $d\geq 1$
		being the spatial dimension. In
		\cite{BethuelsenBirknerDepperschmidtSchluter} an annealed local
		central limit theorem was proven for all $d\geq 1$ and a quenched
		local limit theorem under the assumption $d\geq 3$.
		Here we show that the latter result also holds for all $d \ge 1$.}
	
	\tableofcontents
	
	\section{Introduction and outline of the paper}
	
	This paper concerns the particular model of a directed random walk on
	the backbone of a supercritical oriented percolation cluster. This
	model was first introduced and studied in
	\cite{BirknerCernyDepperschmidtGantert2013}, partly motivated by the study of
	ancestral lineages in evolving populations with locally fluctuating
	population sizes. Since then,
	extensions of both the results and the set-up in
	\cite{BirknerCernyDepperschmidtGantert2013}
	have been addressed by several authors \cite{BethuelsenBirknerDepperschmidtSchluter,
		BirknerCernyDepperschmidtRWDRE2016,
		birkner2024quenched,BirknerGantertSteiber,miller2016random,SteibersPhD2017}.
	
	The aim of the present text is a further extension of the main results of
	\cite{BethuelsenBirknerDepperschmidtSchluter} where certain local
	(central) limit theorems were obtained. In particular, we address the
	extension announced in Remark~1.5 therein by presenting a proof of
	\cite[Theorem 1.3 and 1.4]{BethuelsenBirknerDepperschmidtSchluter}
	that does not pose any restrictions on the dimension. This is inspired
	by work of \cite{peretz2022environment}, who proved an analogous
	extension for certain ballistic random walks in an i.i.d.\ random
	environment, extending the results of \cite{BergerCohenRosenthal2016}.
	
	For more on the motivation for our study and a discussion of related
	literature, we refer to \cite{BethuelsenBirknerDepperschmidtSchluter},
	which (among other things) contains a thorough discussion on how the
	model under study in this paper fits into the more general setting of
	random walk in a (dynamic) random environment. Furthermore, for a more
	detailed discussion of the relevance of our model and results in the context of
	population dynamics, we refer to \cite{BirknerGantert2019}.
	
	The current paper is structured as follows. In
	Section~\ref{sec:model-main-results}, we define the model precisely
	and present the main results. In Section~\ref{sec:thm3.24} we 
	extend the quenched-annealed comparison result \cite[Theorem
	3.24]{SteibersPhD2017} from $d \ge 3$ to all $d \ge 1$. This is the key statement needed in order to
	show that \cite[Theorem~1.3 and 1.4]{BethuelsenBirknerDepperschmidtSchluter}
	can be extended to dimensions $d=1$ and $d=2$. 
	Then, in Section \ref{sec:proofs}, we explain how to apply this
	theorem to obtain our main results, following closely the general
	strategy outlined in \cite[Section~2]{BethuelsenBirknerDepperschmidtSchluter}.

	\section{The model and main results}
	\label{sec:model-main-results}
	
	We start by recalling from \cite[Section
	1.1]{BethuelsenBirknerDepperschmidtSchluter} the precise definition of
	the model and introduce the necessary notation in order to state our
	main results.

	For $d\geq 1$, 
	consider a discrete space-time field $\omega\coloneqq \{\omega(x,n):
	(x,n)\in \bZ^d \times \bZ\}$ of i.i.d.\ Bernoulli random variables
	with parameter $p \in [0,1]$ taking values in $\Omega \coloneqq
	\{0,1\}^{\bZ^d \times \bZ}$, defined on some (large enough)
	probability space equipped with a probability measure $\bP$.  
	For $(x,n) \in \bZ^d \times \bZ$ we  denote by $(x,n)
	\xrightarrow{\omega} \infty$ the event that, for any $m \in \bN$,
	there exists $y \in \bZ^d$ and a space-time sequence $(x_n,n),\dots,
	(x_{n+m},n+m)$ such that $x_n=x$, $x_{n+m}=y$, $\norm{x_k-x_{k-1}} \le 1$ for
	$k=n+1, \dots, n+m$ and $\omega(x_k,k)=1$ for all $k=n,\dots,n+m$
	(here and in the following, $\norm{\cdot}$ denotes the $\sup$-norm).
	
	Let $p_c = p_c(d) \in (0,1)$ 
	be the critical value for this oriented percolation event, i.e.\ the smallest number such that, for all $p>p_c$, it holds that $\bP \big((o,0) \xrightarrow{\omega} \infty
	\big) > 0$, where $o \in \bZ^d$ denotes the origin. 
	We assume throughout this paper that $p>p_c$. 
	The \emph{backbone of a supercritical oriented percolation cluster} as
	alluded to in the title is then defined by
	\begin{align}
		\label{def:cluster}
		\cC \coloneqq \big\{(x,n) \in \bZ^d \times \bZ : (x,n)
		\xrightarrow{\omega} \infty \big\},
	\end{align}
	i.e.\ the set of all space-time sites which are connected to ``time
	$+\infty$'' by a \emph{directed open path}.

	We next define the \emph{(directed) random walk} on $\cC$. That is, given $\omega$, and therefore also the random cluster $\mathcal{C}$,
	this is  the process $(X_{n})_{n\geq 0}$ on $\mathbb{Z}^d$ with
	initial position $X_0=o$ and transition probabilities for $n\ge 0$
	given by
	\begin{align}
		\label{eq:defn_quenched_law}
		\bP(X_{n+1} = y \, \vert \, X_n = x,\omega) =
		\begin{cases}
			\abs{B_1}^{-1}
			& \text{if } (x,n) \notin \mathcal{C}\text{ and } y \in B_1(x);  \\
			\abs{B_1(x)  \cap \mathcal{C}}^{-1}
			& \text{if } (x,n), (y,n+1) \in \mathcal{C} \text{ and } y \in  B_1(x);\\
			0 &\text{otherwise},        
		\end{cases}
	\end{align}
	where,  for a site
	$x \in \bZ^d$ and $l \in \bN$, we write
	$B_l(x) \coloneqq \{ y \in \bZ^d \colon \norm{x-y} \le l \}$ and $B_l\coloneqq B_l(o)$.  
	Thus, the random walk behaves as a simple random walk while not yet on the cluster but 
	stays on the cluster once it is on it.
	We write $P_\omega^{(o,0)}$ for its
	\emph{quenched law}, i.e.\ the conditional law of $\bP$ given
	$\omega$, and $E_\omega^{(o,0)}$ for the corresponding expectation.
	The \emph{annealed (or averaged) law} of $(X_n)$ is denoted by
	$\bP^{(o,0)}$ and its expectation by $\bE^{(o,0)}$. In particular, for
	any $A \in \sigma(X_{n} : n=0,1,\dots)$, we have
	\begin{align}
		\label{eq:15}
		\bP^{(o,0)} (A) = \int P_\omega^{(o,0)}(A)\, d\bP(\omega).
	\end{align}
	Note that whenever we consider the probability of a subset $B\subset\Omega$ 
	we write $\bP(B)$ and there is no superscript since the event $B$ is only 
	concerning the behaviour of $\omega$.
	
	Lastly, we recall the corresponding \emph{environment seen from the
		particle-process}. For this, let $(\sigma_{(y,m)})$ be the
	\emph{space-time shift operator} on $\Omega$ given by $
	\sigma_{(y,m)}\omega(x,n) \coloneqq \omega(x+y,n+m)$.
	Further, denote by $\xi_m(y;\omega) = \xi_m(y)$, where $\xi\coloneqq
	(\xi_n)_{n\in \bZ}$ on $\{0,1\}^{\bZ^d}$ is the process given by 
	\begin{align}
		\label{eq:45}
		\xi_n(x) =\indset{\cC}\bigl((x,n)\bigr), \quad x \in \bZ^d.
	\end{align} 
	The process $\xi$ can be interpreted as the time reversal of the stationary 
	discrete time contact process and therefore is also Markovian and translation invariant.
	Then, a measure $Q$ on $\Omega$ is called \emph{invariant with respect
		to the point of view of the particle} if for every bounded
	continuous function $f:\Omega \to \bR$
	\begin{align}
		\int_\Omega \mathfrak{R} f(\omega)\, dQ(\omega) = \int_\Omega
		f(\omega)\,dQ(\omega),
	\end{align}
	where the transition kernel $\mathfrak{R}$ is given by  
	\begin{align}
		\label{eq:transitionKernel}
		\mathfrak{R}f(\omega) \coloneqq  \sum_{y \in B_1}  g(y;\omega)
		f(\sigma_{(y,1)}\omega),
	\end{align}
	acting on bounded measurable functions $f:\Omega \to \bR$,  and where
	\begin{align}
		\label{eq:19}
		g(y;\omega) \coloneqq
		\ind{\sum\limits_{x \in B_1 } \xi_1(x;\omega) >0,\, \omega(o,0)=1}
		\frac{\xi_1(y;\omega)}{\sum\limits_{x \in B_1} \xi_1(x;\omega)} +
		\ind{\sum\limits_{x \in B_1} \xi_1(x;\omega) =0 \text{ or }
			\omega(o,0)=0 } \frac{1}{3^d}.
	\end{align}
	
	We are now ready to state our first main result. 
	
	\begin{theorem}[{Extension of \cite[Theorem~1.3]{BethuelsenBirknerDepperschmidtSchluter}}]
		\label{thm:main1}
		Let \label{thm:RD for all dimensions} $d\ge 1$ and $p \in (p_c,1]$.
		Then there exists a unique measure $Q$ on $\Omega$ which is
		invariant with respect to the point of view of the particle
		satisfying $Q \ll \bP$ and the concentration property \eqref{eq:concQ}:
		\begin{align}   
			\begin{split} 
				\label{eq:concQ}
				& \text{ there exists $c>0$ so that for every
					$\varepsilon >0$ there is $M_0=M_0(\varepsilon) \in \bN$ and}
				\\ & \bP\Bigl(\Big\lvert\frac{1}{\abs{B_M}}
				\sum_{x\in B_M}\frac{dQ}{d\bP}(\sigma_{(x,0)}\omega)-1
				\Big\rvert>\varepsilon \Bigr) \le M^{-c\log M} \quad \text{ for
					every $M \ge M_0$.}
			\end{split}
		\end{align}
	\end{theorem} 
	
	Our second main result is the following extension of \cite[Theorem
	1.4]{BethuelsenBirknerDepperschmidtSchluter} to also include the cases
	where $d=1$ and $d=2$. For this, letting $Q$ be the measure from
	Theorem \ref{thm:main1}, we denote by $\phi = dQ/d\bP \in L_1(\bP)$
	the Radon–Nikodym derivative of $Q$ with respect to $\bP$.
	
	\begin{theorem}[Quenched local limit theorem, extension of {\cite[Theorem~1.4]{BethuelsenBirknerDepperschmidtSchluter}}]\label{thm:main2}
		Let $d\geq 1$ and $p \in (p_c,1]$. Then, for $\bP$ almost every
		$\omega$, we have
		\begin{align}
			\label{eq:2}
			\lim_{n \rightarrow \infty} \sum_{x\in \bZ^d} \left|
			P_{\omega}^{(o,0)}(X_n = x) - \bP^{(o,0)}(X_n =
			x)\phi(\sigma_{(x,n)} \omega) \right| = 0. 
		\end{align}
	\end{theorem}
	
	\begin{remark}
		\label{rem:ontheclusterornot?}
		In \cite{BirknerCernyDepperschmidtGantert2013} the quenched CLT (and the
		annealed CLT as well as the laws of large numbers) were proven for the
		random walk defined in  \eqref{eq:defn_quenched_law}  using regeneration times (see
		\cite[Eq.~(2.9) in
		Section~2.2]{BirknerCernyDepperschmidtGantert2013} for the precise
		definition of the regeneration sequence). The main point of the construction was
		to \emph{locally} explore the environment and at regeneration times
		(which are not stopping times) find the backbone and determine
		pieces of the path of the random walk between the current and the
		previous regeneration time. In Lemma~2.5 in
		\cite{BirknerCernyDepperschmidtGantert2013} it is shown that the
		increments between the regeneration times have exponential tail
		bounds.
		
		The random walk considered in
		\cite{BirknerCernyDepperschmidtGantert2013} is started on the backbone of the oriented percolation cluster
		whereas we allow in \eqref{eq:defn_quenched_law} the random walk to
		start outside the backbone (the same extension was considered in
		\cite{BethuelsenBirknerDepperschmidtSchluter}). As noted in
		\cite[Remark~2.3]{BirknerCernyDepperschmidtGantert2013} the local
		construction of the regeneration times can be extended to starting
		points outside the cluster. Furthermore, this does not change the
		walk's asymptotic behaviour because when starting away from the
		cluster, the walk will typically hit the cluster very fast, see the
		discussion at the beginning of Section~1.2 in
		\cite{BethuelsenBirknerDepperschmidtSchluter} and also
		\cite[Lemma~B.1]{BethuelsenBirknerDepperschmidtSchluter} for a
		precise quantitative statement.
	\end{remark}
	
	\begin{remark}
		As also alluded to  in \cite[Page 5]{BethuelsenBirknerDepperschmidtSchluter}, it seems conceivable that our main results as stated above can be extended to a larger class of models, e.g.\ those contained within the general framework of \cite{BirknerCernyDepperschmidtRWDRE2016}. One step in this direction was recently obtained in \cite{birkner2024quenched}, where the results of  \cite{BirknerCernyDepperschmidtRWDRE2016} were extended to conclude a quenched central limit theorem.
	\end{remark} 
	
	\section{A first quenched-annealed comparison result}
	\label{sec:thm3.24}
	
	In this section we argue that  \cite[Theorem
	3.24]{SteibersPhD2017}, referred to as Theorem 8.1 in
	\cite{BethuelsenBirknerDepperschmidtSchluter}, holds for $d=1,2$. 
	This is the key technical result that allows us to conclude the main results stated in the previous section.
	This extension is inspired by
	\cite{peretz2022environment}, who extended the results
	of \cite{BergerCohenRosenthal2016} in a similar vein for a related model
	of certain ballistic random walks in random environment. We emphasise that, although we
	apply a similar proof strategy as in
	\cite{BergerCohenRosenthal2016,peretz2022environment}, as  noted in
	\cite{BethuelsenBirknerDepperschmidtSchluter}, there are fundamental
	differences that require new ideas specific to our modelling setting.
	
	Denote by $P_{\omega}^{(y,m)}$ the quenched law with respect to the
	environment $\omega$ of a random walk starting at space-time position
	$z=(y,m) \in \bZ^d \times \bZ$, and by $\bP^{(y,m)}$ the corresponding
	annealed law.  In keeping with the notation from
	\cite{BethuelsenBirknerDepperschmidtSchluter}, we index time for such
	a walk starting from $m$, in particular,
	$P_{\omega}^{(y,m)}(X_n \in \cdot)$ refers to the quenched law of the
	spatial position of that walk after $n-m$ steps for $n \ge m$, etc. 
	Moreover, define for $N \in \bN$,
	\begin{align}
		\label{eq:69}
		\widetilde{\cP}(N) \coloneqq
		\Bigl(\Bigl[-\frac{1}{24}\sqrt{N}\log^3 N,\frac{1}{24}\sqrt{N}\log^3 N \Bigr]^d
		\times \Bigl[0,\frac{1}{3}N\Bigr]\Bigr) \cap (\bZ^d \times \bZ).
	\end{align}
	
	\begin{remark}\label{rem:rescaling}
		We here change the notation compared to \cite{SteibersPhD2017} slightly in that 
		$N$ in \eqref{eq:69} corresponds to $N^2$ in \cite[Eq.~(3.2)]{SteibersPhD2017}. 
		Apart from this slightly different parametrisation, which runs through the whole text, we follow the conventions of  \cite{SteibersPhD2017} closely. 
		For instance, with regard to the definition of $ \widetilde{\cP}(N)$ in \eqref{eq:69}, which corresponds to  \cite[Eq.~(3.2)]{SteibersPhD2017}, we later define an additional set of space-time points that we denote by $\cP(N)$ and that corresponds to \cite[Def.~3.1]{SteibersPhD2017}; see \eqref{eq:defn_P(N)}.
	\end{remark}
	
	For $\theta \in (0,1)$ and $(y,m) \in \widetilde{\cP}(N)$, let $G_4((y,m),N)$
	denote the event that for every box $B \subset \bZ^d$ of side
	length $N^{\theta/2}$ we have
	\begin{align}
		\label{eq:49}
		\big\lvert P_{\omega}^{(y,m)} (X_{N} \in B ) - \bP^{(y,m)}(X_{N}
		\in B) \big\rvert \le N^{-d(1-\theta)/2- \frac{3}{20}\theta}.
	\end{align}
	The corresponding event is defined in the statement of Theorem 3.24 in  \cite{SteibersPhD2017}, see also  \cite[Theorem 8.1]{BethuelsenBirknerDepperschmidtSchluter}. Note that the factor $\frac{3}{20}$ in the exponent  in  \eqref{eq:49} differs from that in \cite{SteibersPhD2017}, who (after taking into account the rescaling discussed in Remark \ref{rem:rescaling}) uses the factor $\frac{1}{6}$ instead. For its application to prove our main results, however, the essential thing is that this factor is strictly positive. 
	Furthermore, set
	\begin{align}
		\label{eq:48}
		G_4(N)\coloneqq \bigcap_{(y,m) \in \widetilde{\cP}(N)} \bigl(G_4((y,m),N) \cup \{\xi_m(y)=0\}\bigr).
	\end{align}
	Although we do not make this explicit in the notation, note that $G_4((y,m),N)$ and $G_4(N)$ both depend on $\theta$. 
	Here, an environment  $\omega \in G_4(N)$ is  considered ``good'' in the sense that for such $\omega$'s the quenched law is
	suitably close to the annealed law on the level of boxes.
	
	The main result of this section is the following theorem.
	
	\begin{theorem}[Extension of  {\cite[Theorem~3.24]{SteibersPhD2017}}]
		Let \label{thm: Theorem 3.24 Ersatz} $d\ge 1$ and $\theta \in
		(0,1)$. There exist positive constants $c$ and $C$, such that for
		all $(y,m) \in \widetilde{\cP}(N)$  with $N \in \bN$  we have
		\begin{align}
			\bP\bigl(G_4((y,m),N)  \cup \{\xi_m(y)=0\} \bigr) \ge 1 - C N^{-c\log N}
		\end{align}
		and
		\begin{align}
			\label{eq: lower bound for G4(N)}
			\bP\bigl(G_4(N) \bigr) \ge 1- CN^{-c\log N}.
		\end{align}
	\end{theorem}

	Theorem~\ref{thm: Theorem 3.24 Ersatz} is the extension of {\cite[Theorem~3.24]{SteibersPhD2017}} from $d\geq 3$ to $d\geq 1$. Equipped with this statement, the proofs of our main results follow similarly to the arguments of  \cite{BethuelsenBirknerDepperschmidtSchluter} after minor modifications. See the following section for more details.

	The proof of Theorem~\ref{thm: Theorem 3.24 Ersatz} hinges on the
	following lemma. To state this, we define another class of good events. In particular, for $N \in \bN$, $z\in \widetilde{\cP}(N)$ 
	and $\theta \in (0,1)$, let $G_1(z,\theta,N)$ be the event that, for
	every  $M \in [\frac{2}{5}N, N]$ and every box $B \subset
	[-\sqrt{N}\log^3(N),\sqrt{N}\log^3(N)]^d$ of side length $N^{\theta/2}$, we have
	\begin{align}
		\label{eq:ann-que.thetabox}
		\abs{P^z_\omega(X_M\in B) - \bP^z(X_M\in B)}\le N^{d(\theta-1)/2}.
	\end{align}
	Here one should note that by the annealed local CLT
	\cite[Theorem~1.1]{BethuelsenBirknerDepperschmidtSchluter}, we have  $\bP^z(X_M\in B) = O((N^{\theta/2})^d \cdot N^{-d/2})$  for any $\theta \in (0,1)$. 
	Thus \eqref{eq:ann-que.thetabox} is in itself not a particularly tight
	bound, but rather a statement that (on the level of boxes of diameter
	$N^{\theta/2}$) the quenched law is of comparable smoothness as the
	annealed law. Nonetheless, this property is an important ingredient
	in the proof of Theorem~\ref{thm: Theorem 3.24 Ersatz}, which
	establishes much tighter bounds of the form \eqref{eq:49} on the
	difference between the quenched and the annealed law on boxes.
	See also the discussion in \cite[Remark~3.12]{SteibersPhD2017}.
	
	\begin{lemma}[Extension of Proposition~3.11 in \cite{SteibersPhD2017}]
		Let $d\ge1$ and recall the event $G_1(z,\theta,N)$ defined with \eqref{eq:ann-que.thetabox}. \label{lem:Prop 3.11 Ersatz} 
		Then, for every $\frac{2d+1}{2d+2}<\theta\le 1$ there exist
		constants $C,c>0$, such that
		\begin{align}
			\label{eq:7}
			\bP(G_1(z,\theta,N)) \ge 1-CN^{-c\log N}
		\end{align}
		for all  $N \in \bN$ and  $z\in \widetilde{P}(N)$, and therefore also
		\begin{align}
			\label{eq:8}
			\bP\bigg(\bigcap_{z\in\widetilde{P}(N)} G_1(z,\theta,N) \bigg) \ge 1- CN^{-c\log N}.
		\end{align}
	\end{lemma}

	Lemma \ref{lem:Prop 3.11 Ersatz} extends {\cite[Proposition~3.11]{SteibersPhD2017}} to  $d=1$ and $d=2$. 
	For this, in the statement of  Lemma \ref{lem:Prop 3.11 Ersatz},  we make a stronger assumption on the value of $\theta$ compared to {\cite[Proposition~3.11]{SteibersPhD2017}}, where the equivalent statement for $d\geq 3$ was shown for every $\frac{d}{d+1}<\theta\le 1$. 
	Our requirement on $\theta$ in Lemma \ref{lem:Prop 3.11 Ersatz} stems from the fact that in our adaptation of {\cite[Corollary~3.20]{SteibersPhD2017}},  we use a much more conservative bound on the number of times that two independent random walkers in the same environment  are allowed to meet before time $N$, namely $N^{1/2+\varepsilon}$ instead of $\log^{12}(N)$; see \eqref{eq:replacementLemma} in Lemma~\ref{lem:replacement_cor_3.20}.  
	We want to highlight here that it is exactly this change that allows us to prove Lemma~\ref{lem:Prop 3.11 Ersatz} for $d=1,2$, albeit under  the slightly more restrictive assumption 
	that $\theta \geq \frac{2d+1}{2d+2}$ that we apply  precisely at Equation \eqref{eq:lower_bound_theta} in the proof.

	The remainder of this section is organised as follows. 
	In the next subsection we first collect a couple of statements from
	\cite{BethuelsenBirknerDepperschmidtSchluter} and \cite{SteibersPhD2017}. 
	Fundamental to the proof of Lemma \ref{lem:Prop 3.11 Ersatz}  is the extension of {\cite[Corollary~3.20]{SteibersPhD2017}} to $d=1$ and $d=2$, detailed in Subsection \ref{sec:lem:replacement_cor_3.20}. 
	Then, in Subsection \ref{sec:proof_of_lemma_3.3}, we present the proof of  Lemma \ref{lem:Prop 3.11 Ersatz}.
	The proof of Theorem \ref{thm: Theorem 3.24 Ersatz} is summarized in the last subsection. 
	
	\begin{remark}
		In the following proofs, $C$ and $c$ denote unspecified constants whose
		precise value is not relevant and may even change from line to line.
	\end{remark}

	\subsection{Preliminary annealed estimates}
	
	The following statements from
	\cite{BethuelsenBirknerDepperschmidtSchluter} and
	\cite{SteibersPhD2017} were formally stated therein for $d\geq 3$, but
	their proofs are all based on estimates stemming from the annealed CLT
	of \cite[Theorem~1.1]{BirknerCernyDepperschmidtGantert2013} obtained
	for $d\geq 1$, and hence hold in all dimensions. 
	
	\begin{lemma}[Extension of Lemma 3.1 in \cite{BethuelsenBirknerDepperschmidtSchluter}]
		For any \label{lem:annealed_derivative_estimates} $d \ge 1$, 
		there is a constant $C=C(d)<\infty$, such that, for any 
		$j=1,\dots,d$, $x,y\in \bZ^d$, $m,n \in \bN$, 
		\begin{align}
			\label{eq:1}
			\abs{\bP^{(y,m)}(X_{n+m} = x) -\bP^{(y+e_j,m)}(X_{n+m} = x)}
			& \le Cn^{-(d+1)/2},\\
			\label{eq:3}
			\abs{\bP^{(y,m)}(X_{n+m} = x) -\bP^{(y,m+1)}(X_{n+m} = x)}
			& \le Cn^{-(d+1)/2},\\
			\label{eq:4}
			\abs{\bP^{(y,m)}(X_{n+m} = x) -\bP^{(y,m)}(X_{n+m} = x+e_j)}
			& \le Cn^{-(d+1)/2},\\
			\label{eq:5}
			\abs{\bP^{(y,m)}(X_{n+m} = x) -\bP^{(y,m)}(X_{n-1+m} = x)}
			&\le Cn^{-(d+1)/2},
		\end{align}
		where  $e_j$ denotes the $j$-th (canonical) unit vector.
	\end{lemma}
	
	\begin{lemma}[Extension of Lemma 3.2 in \cite{BethuelsenBirknerDepperschmidtSchluter}]
		\label{lem:additional_annealed_estimate} Let $d\ge 1$ and $\varepsilon>0$. For
		$n\in\bN$ large enough and every partition $\Pi^{(\varepsilon)}_n$
		of $\bZ^d$ into boxes of side length $\lfloor n^\varepsilon
		\rfloor$, we have
		\begin{align}
			\label{eq:6}
			\sum_{B \in \Pi^{(\varepsilon)}_n}\sum_{x\in B}
			\max_{y\in B}\bigl[\bP^{(o,0)}(X_n=y)-\bP^{(o,0)}(X_n=x)\bigr]
			\le Cn^{-\frac{1}{2}+ 3d\varepsilon}, 
		\end{align}
		with a constant $C=C(d,\varepsilon)<\infty$.
	\end{lemma}

	\begin{lemma}[Lemma~3.6 from \cite{SteibersPhD2017}]\label{lem:lem3.6SS}
		Let $d \ge 1$, $(y,m) \in \widetilde{\cP}(N)$ and $n \in [N/2,N]$.
		There exist positive constants $c, C>0$, depending on $d$ but not
		on $(y,m)$ and $n$, so that, for all $N \in \bN$, 
		\begin{align}
			\label{eq:9}
			\mathbb P^{(y,m)} (\norm{X_n - y} \ge \sqrt{n} \log^3(N)) \le C
			N^{-c \log (N)}. 
		\end{align}
		Furthermore, let $H((y,m),N)$ be the event that
		\begin{align}
			\label{eq:11}
			P^{(y,m)}_\omega (\norm{X_n - y} \ge \sqrt{n} \log^3(N)) \le C
			N^{-c \log (N)/2},  
		\end{align}
		for all $N/2\le n\le N$. Then 
		\begin{align}
			\label{eq:10}
			\mathbbm P\bigl(H((y,m),N) \bigr) \ge 1 - C N^{-c \log
				(N)/2}. 
		\end{align}
	\end{lemma}
	
	The key property applied in order to prove the above statements is that, as we also noted in Remark \ref{rem:ontheclusterornot?}, the random walk $(X_n)$ possesses a regeneration
	construction.  Thus, under the annealed law, it is literally a centred
	random walk when observed along the regeneration times
	where the inter-regeneration increments in fact have exponential
	tails. 
	Particularly, when observed along regeneration times,
	estimates analogous to \eqref{eq:1}--\eqref{eq:5} above are standard
	consequences of the classical local CLT.
	For the convenience of the reader we provide more proof details for the above lemmas
	in Appendices\ \ref{sec:pf_lem:annealed_derivative_estimates}--\ref{sec:pf_lem3.6SS}.
	
	\subsection{Bounds on the intersection of two random walks}\label{sec:lem:replacement_cor_3.20}

	We now present the aforementioned extension of \cite[Corollary 3.20]{SteibersPhD2017} to $d=1$ and $d=2$, which is the essential change in our approach compared to \cite{SteibersPhD2017}. For this, we denote by  $0=T^\simu_0<T^\simu_1<\dots$, the  simultaneous regeneration times introduced in \cite[Equation (3.7)]{BirknerCernyDepperschmidtGantert2013}  for two independent random walks, $X$ and $X'$, evolving in the same environment and initiated at $X_0=y$ and $X'_0=y'$, respectively. 
	The precise definition of these regeneration times is not important here and we will in the following only use (beside the regeneration property) the fact that 
	\begin{align}
		\label{eq:regen_tail}
		\bP(T^\simu_i-T^\simu_{i-1}>k ) \le C\exp(-ck), \quad k > 0,
	\end{align}
	as concluded in \cite[Lemma 3.1]{BirknerCernyDepperschmidtGantert2013}. Note that the probability
	in \eqref{eq:regen_tail} in actuality depends on the two starting positions $X_0=y$ and $X'_0=y'$. However, since
	the upper bound is uniform for any choice of $y$ and $y'$, wo omit denoting them at the probability measure $\bP$.
	Further, we consider the event
	\begin{equation}
		\label{eq:defn_RN}
		R^{\simu}_N \coloneqq \{ T^\simu_k-T^\simu_{k-1}\le \log^2(N) \text{ for all } k\le N \}.
	\end{equation}
	We are now ready to state the main result of this subsection.
	
	\begin{lemma}[Replacement for Corollary~3.20 from \cite{SteibersPhD2017} for $d=1,2$]
		\label{lem:replacement_cor_3.20}
		For any $d\ge 1$ and $\varepsilon>0$ 
		there are constants $C,c>0$ such that, irrespective of the initial positions $y,y'\in \bZ^d$, 
		\begin{equation}
			\label{eq:replacementLemma}
			\bP\bigg(\Big\{ \sum_{i=1}^N \ind{\abs{X_i-X'_i}<4\log^2(N)} > N^{1/2+\varepsilon}\Big\} \cap R^\simu_N \bigg) \le CN^{-c\log N}
		\end{equation}
		for all $N \in \bN$.
	\end{lemma}

	\begin{proof}[Proof of Lemma~\ref{lem:replacement_cor_3.20}]
		
		The statement for $d\geq 3$ follows immediately from \cite[Corollary 3.20]{SteibersPhD2017} since $N^{1/2+\varepsilon} > \log^{12}(N)$. (Alternatively, the present proof could easily be adapted to the case $d \ge 3$ as well.)
		
		The proof for $d=2$ builds on the proof of \cite[Lemma~3.9]{BirknerCernyDepperschmidtGantert2013}. 
		Particularly, as in \cite[Equation (3.39)]{BirknerCernyDepperschmidtGantert2013}, we define stopping times $\mathcal{R}_i$, $\mathcal{D}_i$ and $\cU$ by $\mathcal{R}_0=0$, and for $i\geq 1$,
		\begin{align}
			\label{eq:defn_Di_Ri_U}
			\begin{split}
				\mathcal{D}_i &\coloneqq \inf\{ k\ge \mathcal{R}_{i-1} \colon \norm{\widehat{X}_k - \wh{X}^{'}_k } \ge N^{b_1} \},\\
				\mathcal{R}_i &\coloneqq \inf\{ k\ge \mathcal{D}_i \colon \norm{\wh{X}_k -\wh{X}^{'}_k } \le \log^4(N) \},\\
				\mathcal{U} &\coloneqq \inf\{ k\ge 0 \colon \norm{\wh{X}_k -\wh{X}^{'}_k } \ge N \},
			\end{split}
		\end{align}
		where 
		$\wh{X}_k \coloneqq X_{T^\simu_k}$ and $\wh{X}_k^{'} \coloneqq X'_{T^\simu_k}$ are the two walks observed along their simultaneous regeneration times, and 
		$b_1 \in (0,1/2)$ is such that the statement of \cite[Lemma~3.8]{BirknerCernyDepperschmidtGantert2013} holds (Comparing with \cite[p.~28]{BirknerCernyDepperschmidtGantert2013} we observe that $\mathcal{R}_i$ is defined
		there slightly differently, using $K \log N$ instead of $\log^4(N)$. However, this does not affect our arguments). Then, the argumentation leading up to \cite[Inequality (3.41)]{BirknerCernyDepperschmidtGantert2013} also gives that
		\begin{equation}\label{eq:help1+2}
			\bP \bigg(\sum_{i=1}^{J} \mathcal{D}_i-\mathcal{R}_{i-1} \ge N^{1/2+\varepsilon}\log^{-2}(N) \bigg) \leq C\exp(-cN^{c'}),
		\end{equation}
		where $J$ is the unique integer such that $\mathcal{D}_J \le \mathcal{U}\le \mathcal{R}_J$. 
		To briefly recall the idea to prove this note that $J$ is a geometric random variable with parameter converging to $b_1$ as $N\to \infty$ (this can be shown using Lemma~3.6 from \cite{BirknerCernyDepperschmidtGantert2013}, see also the proof of Lemma~3.9 in \cite{BirknerCernyDepperschmidtGantert2013}) which gives an upper bound on the number of summands, i.e. there exist constants $C,c>0$ such that for any $\delta>0$ we have
		\begin{equation*}
			\bP(J>N^\delta) \le C\exp(-cN^\delta).
		\end{equation*}	
		Then one can use the separation lemma, see Lemma~3.8 in \cite{BirknerCernyDepperschmidtGantert2013}, to get the upper bound
		\begin{equation*}
			\bP\big(\mathcal{D}_i-\mathcal{R}_{i-1} >N^{b_2}  \big) \le e^{-b_3N^{b_4}},
		\end{equation*}
		where $b_2,b_3$ and $b_4$ are positive constants from \cite[Lemma~3.8]{BirknerCernyDepperschmidtGantert2013}. The relation between $b_1$ and $b_2$ is crystallized in Step~6 of the proof of Lemma~3.8 in \cite{BirknerCernyDepperschmidtGantert2013}, where the authors need that $(2b_4+3b_1)\vee (3b_4) < b_2$. Also note that the bound in this separation lemma is universal in the starting positions which is also the reason that \eqref{eq:replacementLemma}   does not depend on   
		the starting positions. To prove \eqref{eq:help1+2} it is thus sufficient to have $b_2<1/2$, which we can obtain by choosing $b_1$ and $b_4$ small enough.
		
		Now, on the event $R^\simu_N$, the time between two consecutive regeneration times is smaller than $\log^2(N)$ up until the $N$-th regeneration time.  
		Further, the number of steps in the original pair of random walks $(X,X')$ at which they are closer than $\log^2(N)$ until time $N$ can be bounded from above by observing the pair $(\wh{X},\wh{X}')$. Particularly, we count the times at which these are closer than $\log^4(N)$ until time $J$ and multiply that number by $\log^2(N)$. Then, since $\log^4(N)-\log^2(N)>4\log^2(N)$, on $R^\simu_N$,  we obtain that
		\begin{equation*}
			\sum_{i=1}^{N} \ind{\abs{X_i-X'_i}<4\log^2(N)} \le \log^2(N)\sum_{i=1}^J \big(\mathcal{D}_i-\mathcal{R}_{i-1}\big).
		\end{equation*}
		It thus follows that
		\begin{align*}
			\bP&\bigg( \Big\{\sum_{i=1}^N \ind{\abs{X_i-X'_i}<4\log^2(N)} > N^{1/2+\varepsilon}\Big\} \cap R^\simu_N \bigg) \\
			&\le \bP\bigg( \Big\{\log^2(N)\sum_{i=1}^J \mathcal{D}_i-\mathcal{R}_{i-1} \ge N^{1/2+\varepsilon} \Big\} \cap R^\simu_N \bigg)\\
			&\le \bP\bigg( \sum_{i=1}^J \mathcal{D}_i-\mathcal{R}_{i-1} \ge N^{1/2+\varepsilon}\log^{-2}(N)  \bigg)\\
			&\le C\exp(-cN^{c'}),
		\end{align*}
		where we have applied \eqref{eq:help1+2} in the last inequality. This concludes the proof for $d=2$.

		For the case $d=1$ we follow the analysis of \cite[Lemma~3.14]{BirknerCernyDepperschmidtGantert2013}.
		Recall the definitions of $\cR_i$ and $\cD_i$ from \eqref{eq:defn_Di_Ri_U} and that $\cR_0=0$.
		Note that $b_1$ is a parameter that has to be tuned correctly, depending on the $\varepsilon>0$ in the statement of Lemma~\ref{lem:replacement_cor_3.20}. 
		In addition, let $I_n\coloneqq \max\{ i\colon \mathcal{R}_{i}\le n \}$ be the number of ``black boxes'' up to time $n$, i.e.\ the number of regeneration intervals in which we cannot guarantee a minimum distance of $\log^2(N)$ between the two walks, even if the event $R^\simu_N$ occurs.
		We recall a property of $I_n$, see \cite[Equation 3.95]{BirknerCernyDepperschmidtGantert2013}, namely 
		\begin{equation}
			\bP(I_n\ge k) \le \exp(-ck^2/n)
		\end{equation}
		for $1\le k\le n$ and some constant $c>0$. On $R^\simu_N$, and by the definition of $I_n$, by the same arguments as for the case $d=2$, we obtain that  
		\begin{equation*}
			\sum_{i=1}^{N} \ind{\abs{X_i-X'_i}<4\log^2(N)} \le \log^2(N)\sum_{i=1}^{I_{N}+1} (\mathcal{D}_{i}-\mathcal{R}_{i-1}).
		\end{equation*}
		Therefore, it follows that 
		\begin{align*}
			\bP&\bigg( \Big\{\sum_{i=1}^N \ind{\abs{X_i-X'_i}<4\log^2(N)} > N^{1/2+\varepsilon}\Big\} \cap R^\simu_N \bigg) \\
			&\le \bP\bigg( \Big\{ \log^2(N)\sum_{i=1}^{I_{N}+1} \mathcal{D}_{i}-\mathcal{R}_{i-1} > N^{1/2+\varepsilon}\Big\} \cap R^\simu_N \bigg)\\
			&\le \bP(I_{N}\ge N^{1/2+\varepsilon/2}) +\bP\bigg( \log^2(N)\sum_{i=1}^{I_{N}+1} \mathcal{D}_{i}-\mathcal{R}_{i-1} > N^{1/2+\varepsilon} \,,\, I_{N}+1 < N^{1/2+\varepsilon/2}\bigg)\\
			&\le\exp(-cN^{\varepsilon}) + \sum_{i=1}^{N^{1/2+\varepsilon/2}} \bP\bigg( \mathcal{D}_{i}-\mathcal{R}_{i-1} > N^{\varepsilon/2}/\log^2 (N) \bigg).
		\end{align*}
		Then, by choosing $b_1 < \varepsilon/2$ in the definition of $\mathcal{D}_i$, the last term can be bounded using \cite[Lemma~3.15]{BirknerCernyDepperschmidtGantert2013} 
		which concludes the proof for $d=1$.
		Even though it is not stated there explicitly, the exponent $b_2 > 0$ in \cite[Eq.~(3.88) in Lemma~3.15]{BirknerCernyDepperschmidtGantert2013}, which is analogous to our $b_1$ here, can be taken arbitrarily small. Arguably, the argument in \cite{BirknerCernyDepperschmidtGantert2013} is rather sketchy. We refer to \cite[Lemma~3.31 and its proof in Section~A.6]{birkner2024quenched} for a detailed discussion of the crucial separation lemma in case $d=1$. Note that the model setup in \cite{birkner2024quenched} is a little different from \cite{BirknerCernyDepperschmidtGantert2013} but the arguments and tools used in this specific proof are completely parallel.
	\end{proof}
	
	\subsection{Proof of Lemma \ref{lem:Prop 3.11 Ersatz}}
	\label{sec:proof_of_lemma_3.3}
	
	Most of the arguments used in the proof of \cite[Proposition 3.11]{SteibersPhD2017} hold for $d\ge 1$ without any change. We introduce the notation and recall some estimates from there that we need for the proof of Lemma \ref{lem:Prop 3.11 Ersatz}. For this, we first enumerate the elements of
	\begin{equation}
		\label{eq:defn_P(N)}
		\mathcal{P}(N) \coloneqq \big(\bZ^d\times\bZ \big)\cap \big( [ -\sqrt{N}\log^3(N),\sqrt{N}\log^3(N) ]^d\times[0,N+\log^3 (N)] \big)
	\end{equation} 
	by ordering them increasing in time and then lexicographically in space.
	Now, fix $\theta > \frac{2d+1}{2d+2}$ and let $\theta' \in (\frac{2d+1}{2d+2}, \theta)$. Further, for $z\in \mathcal{P}(N)$ and $v\in\bZ^d$, define for $M \in [\frac{2}{5}N, N]$
	\begin{equation}
		\label{eq:defn_martingale}
		M_k \coloneqq \bE \Big[ P^z_\omega (X_{M+N^{\theta'}}=v)\,\vert\, \mathcal{F}_k \Big], \qquad k=0,\dots, |\cP(N)|,
	\end{equation}
	where $\mathcal{F}_k \coloneqq \sigma(\omega(z_1),\dots,\omega(z_k))$ is the $\sigma$-algebra generated by the  $\omega$'s on the first $k$ elements of $\mathcal{P}(N)$. Note that we use the variable $M$ in the definition of $M_k$. Since $M$ will always play the same role and is often a fixed parameter in proofs, we omit specifically denoting this dependence on the left hand side of \eqref{eq:defn_martingale}. Moreover, by definition, $M_k=M_k(\omega)$ is a random variable, i.e.\ a (real-valued) function on $\Omega$, which in fact only depends on the value of $\omega$ at the space-time coordinates $z_1,\dots,z_k$ (see the proofs of Lemmas~3.15 and 3.18 in \cite{SteibersPhD2017} for a more explicit description of this function). When we want to make this fact notationally explicit in the computations below, we will write $\bE[P^z_{\widetilde{\omega}}(X_{M+N^{\theta'}}=v)\vert \mathcal{F}_k ](\omega)$, 
	where the $P^z_{\widetilde{\omega}}$ inside the conditional expectation signals that we are averaging over all coordinates of $\mathcal{P}(N) \setminus  \{z_1,\dots,z_k\}$.
	
	Throughout this subsection $z$ is a fixed starting position in space-time, whereas $v$ is a fixed position in space. However, note that all estimates obtained in the following lemmas hold with a uniform bound for all $v$, see Lemma~\ref{lem:McDiarmidApplied}. We make it clear whenever we restrict $v$ to a certain region.
	Moreover, consider the sequence $(U_k)_{k=0}^{|\mathcal{P}(N)|}$ given by   
	\begin{equation}
		\label{eq:defn_Uk}
		U_k \coloneqq \text{esssup}\Big( \Abs{M_{k-1}-
			M_{k}}\,\Big\vert \mathcal{F}_{k-1} \Big)(\widetilde{\omega}).
	\end{equation}
	Thus, $U_k$ controls how much the (partially) quenched probability changes if we expose one additional site of $\omega$, namely $z_k$. 
	The idea of the proof of Lemma \ref{lem:Prop 3.11 Ersatz} is now to telescope over the space-time points of $\cP(N)$ using the above enumeration and applying the Martingale inequality of \cite[Thm.~3.14]{McDiarmid1998}:
	\begin{equation}\label{ineq:McDiarmid}
		\bP\big( \abs{M_n - M_0}\ge \alpha, U\le c \big) \le 2\exp\Big(-\frac{\alpha^2}{2c} \Big), \quad \text{ where } U \coloneqq \sum_{k=1}^{|\cP(N)|} U_k^2.
	\end{equation}
	In order to apply this inequality we first need a bound on the terms $U_k$ that we present next. Here, with $z_k=(y_k,m_k)$ denoting the $k$'th element in the enumeration of $\cP(N)$, we let  
	\begin{equation}
		\label{eq:Delta_k}
		\Delta_k=\Delta_k(z_k,N) \coloneqq \{ (y,m)\in\mathcal{P}(N) \colon \norm{y-y_k}\le m_k-m \le \log^2 (N) \}
	\end{equation}
	and write $\{\Delta_k \text{ visited} \}$ to denote the event that the random walk visits the space-time region $\Delta_k$. (See  \cite[Page 55, Figure 3.2]{SteibersPhD2017} for a graphical representation of the set $\Delta_k$, there named $\textbf{M}_k$)

	\begin{lemma}\label{lem:BoundingU_k}
		For $z\in \mathcal{P}(N)$ there is an event $\widetilde{G}_1(z,N)$ and constants $C,c <\infty$, such that for all  $\widetilde{\omega}\in\widetilde{G}_1(z,N)$  
		\begin{align}
			\label{eq:BoundingU_k}
			U_k(\widetilde{\omega}) \le C\log^{6d+9}(N) P^z_{\widetilde{\omega}}(\Delta_k \text{ visited})N^{-(d+1)\theta'\!/2} + CN^{-c\log N} \text{ for all }k \le \abs{\mathcal{P}(N)}
		\end{align}
		and satisfying  $\bP(\widetilde{G}_1(z,N))\ge 1-CN^{-c\log N}$, where $\theta'$ is from \eqref{eq:defn_martingale}.
	\end{lemma}

	Lemma \ref{lem:BoundingU_k} can be seen as a direct consequence of Lemma~3.15 and Lemma~3.18 in \cite{SteibersPhD2017}. See the end of this subsection for a description of its proof.

	As in the previous subsection, we let $X'$ be an additional random walk, evolving independently of $X$ but in the same environment and initiated at the same location. 
	Consider the event 
	\begin{equation}
		\label{eq:defn_W(N)}
		W(N) \coloneqq  \Big\{ \sum_{i=1}^N \ind{\norm{X_i-X'_i}<4\log^2(N)} \le N^{1/2+\varepsilon}\Big\} \cap R^\simu_N.
	\end{equation}
	By a slight abuse of notation, in the following we also write $P_{\omega}^z$ for the product measure $P_\omega^z\otimes P_\omega^{z'}$ whenever we consider both walkers $X$ and $X'$ started at the space-time locations $z$ and $z'$, respectively.
	
	Moreover, for $M \in [\frac{2}{5}N,N]$, 
	we define $k_0$ to be the value in the enumeration of $\cP(N)$  
	such that $m_{k_0} = M+\log^2 N$ and $(y_k,m_k)=z_k \prec z_{k_0}=(y_{k_0},m_{k_0})$ for all $z_k$ such that $m_k \le M+\log^2 N$. Utilising Lemma \ref{lem:BoundingU_k} and the inequality in \eqref{ineq:McDiarmid}, we obtain the following. 
	
	\begin{lemma}\label{lem:McDiarmidApplied}
		For all $\theta'\in (\tfrac{2d+1}{2d+2},\theta)$ there exist constants $c,C>0$ such that, uniformly in $v\in\bZ^d$ and $M \in [\frac{2}{5}N,N]$,
		\begin{equation*}
			\bP\Big( \abs{\bE[P^z_{\widetilde{\omega}}(X_{M+N^{\theta'}}=v)\vert \mathcal{F}_{k_0} ](\omega)- \bP^z(X_{M+N^{\theta'}}=v) }> \frac{1}{2}N^{-d/2} \Big)
			\le C N^{-c\log N}.
		\end{equation*}
	\end{lemma} 
	
	\begin{proof}
		Let  $\omega \in \widetilde{G}_1(z,N)\cap W(N)$ with $\widetilde{G}_1(z,N)$ from Lemma~\ref{lem:BoundingU_k}, $W(N)$ as in \eqref{eq:defn_W(N)} and denote by $\ell(N)$ a slowly varying function which we use to abbreviate all the $\log$ factors in the display below. Then, by  Lemma \ref{lem:BoundingU_k}, we have that 
		\begin{align}
			\notag
			\sum_{k=1}^{k_0} U_k^2 &\le CN^{-c\log N} + C\ell(N) N^{-(d+1)\theta'} \sum_{k=1}^{k_0} (P_\omega^z (\Delta_k \text{ visited}))^2\\
			\notag
			&\le CN^{-c\log N} + C\ell(N)N^{-(d+1)\theta'}\sum_{j=1}^{M+\log^2 N} P^z_\omega \Big( \norm{X_j-X'_j}< 4\log^2 (N) \Big)\\
			\label{eq:12}
			&\le CN^{-c\log N} + C\ell(N)N^{-(d+1)\theta'}N^{1/2+\varepsilon},
		\end{align}
		where the second inequality follows by noticing that 
		\begin{equation}
			\big(P^z_\omega (\Delta_k \text{ visited})\big)^2 = P^z_{\omega}(\Delta_k \text{ visited by both }X \text{ and } X')
		\end{equation}
		together with the fact that 
		if both random walks hit $\Delta_k$ they are necessarily at distance at most $4\log^2 (N)$ at time $m_k$, cf.\ \eqref{eq:Delta_k}. Note also that for $(y_k,m_k)$ and $(y_{k'},m_{k'})$ with $k \neq k'$ and $m_k=m_{k'}$, we have $\Delta_k \cap \Delta_{k'} = \emptyset$ whenever $\norm{y_k-y_{k'}}>2\log^2 N$.
		The last inequality in \eqref{eq:12} holds true by the definition of $W(N)$.  
		Therefore, using McDiarmid's inequality from \eqref{ineq:McDiarmid}, we get
		\begin{multline*}
			\bP\Big( \abs{\bE[P^z_{\widetilde{\omega}}(X_{M+N^{\theta'}}=v)\vert \mathcal{F}_{k_0} ](\omega)- \bP^z(X_{M+N^{\theta'}}=v) }> \frac{1}{2}N^{-d/2} \Big)\\
			\le 2\exp \Big( -\frac{CN^{-d}}{\ell(N)N^{-(d+1)\theta'+1/2+\varepsilon}}\Big) + \bP\big(\widetilde{G}_1(z,N)^\compl\cup W(N)^\compl\big).
		\end{multline*}
		Since, for some $\varepsilon>0$, it holds that 
		\begin{equation}
			\label{eq:lower_bound_theta}
			\theta'>\frac{d+\frac{1}{2}+\varepsilon}{d+1} = \frac{2d+1+2\varepsilon}{2d+2}
		\end{equation}
		and hence $-d+(d+1)\theta'-\frac{1}{2}-\varepsilon >0,$ this, together with Lemma~\ref{lem:BoundingU_k} and Lemma~\ref{lem:replacement_cor_3.20}, completes the proof.
	\end{proof}
	
	Finally, combining the above estimates, we are now ready to present the proof of Lemma \ref{lem:Prop 3.11 Ersatz}.

	\begin{proof}[Proof of Lemma \ref{lem:Prop 3.11 Ersatz}]
		Recall $k_0$ from above Lemma~\ref{lem:McDiarmidApplied} and let $\widetilde{G}_2(z,N)$ be the event that
		\begin{equation*}
			\abs{\bE[P^z_{\widetilde{\omega}}(X_{M+N^{\theta'}}=v) \,\vert\, \mathcal{F}_{k_0}](\omega)-\bP^z(X_{M+N^{\theta'}}=v)} \le \frac{1}{2}N^{-d/2}
		\end{equation*}
		for every $M\in [\frac{2}{5}N,N]$ and every $v \in B_{\sqrt{N}\log^3(N)}$,
		where $\theta' \in (\frac{2d+1}{2d+2},\theta)$ as in \eqref{eq:defn_martingale}.
		Then, since the number of $v$'s we consider here is bounded by $(2\sqrt{N}\log^{3d}(N)+1)^d$, we obtain from Lemma~\ref{lem:McDiarmidApplied} that $\bP(\widetilde{G}_2(z,N))\ge 1-CN^{-c\log N}$. 
		
		Now, fix a $\omega \in \widetilde{G}_2(z,N)$ and  let $B_x^{(0)} = B_{N^{\theta/2}}(x)$ with $x\in \bZ^d$.
		Then, as in \eqref{eq:ann-que.thetabox}, we want to estimate $\abs{P^z_\omega(X_M \in B_x^{(0)})-\bP^z(X_M \in B_x^{(0)})}.$  
		For this, let $B^{(1)}_x = B_{\frac{9}{10}N^{\theta/2}}(x)$ and $B^{(2)}_x=B_{\frac{11}{10}N^{\theta/2}}(x)$ and note that,  
		by Lemma \ref{lem:lem3.6SS} and since $\theta'<\theta$, it follows that
		\begin{align*}
			\bP^z(X_{M+N^{\theta'}}\in B^{(1)}_x) &\leq  \bP^z(X_{M}\in B_x^{(0)}) + \bP^z(X_{M+N^{\theta'}}\in B^{(1)}_x, X_{M}\notin B_x^{(0)})\\
			&\le \bP^z(X_{M}\in B_x^{(0)}) + \bP^z(\norm{X_{M+N^{\theta'}}-X_M}>\frac{1}{10}N^{\theta/2})\\
			&= \bP^z(X_{M}\in B_x^{(0)}) + \bP^{(o,0)}(\norm{X_{N^{\theta'}}}>\frac{1}{10}N^{\theta/2})\\
			&\le \bP^z(X_{M}\in B_x^{(0)}) + \bP^{(o,0)}(\norm{X_{N^{\theta'}}}>\sqrt{N^{\theta'}}\log^3(N))\\
			&\le  \bP^z(X_{M}\in B_x^{(0)}) + CN^{-c\log N},
		\end{align*}
		where we used translation invariance in the third line, and therefore
		\begin{align}
			\label{eq:14}
			\bP^z(X_{M+N^{\theta'}}\in B^{(1)}_x) &< \bP^z(X_{M}\in B_x^{(0)}) + CN^{-c\log N},\\
			\shortintertext{and analogously}
			\label{eq:16}
			\bP^z(X_{M+N^{\theta'}}\in B^{(2)}_x) &>\bP^z(X_{M}\in B_x^{(0)}) - CN^{-c\log N}.
		\end{align}
		Additionally, there exists a set $\widetilde{G}_3(z,N)$ such that for all $\widetilde{\omega}\in \widetilde{G}_3(z,N)$
		\begin{align}
			\label{eq:17}
			\bE\Big[ P^z_\omega(X_{M+N^{\theta'}}\in B^{(1)}_x)\,\vert\, \mathcal{F}_{k_0} \Big] (\widetilde{\omega}) &< P^z_{\widetilde{\omega}}(X_M \in B_x^{(0)}) + CN^{-c\log N},\\
			\label{eq:18}
			\bE\Big[ P^z_\omega(X_{M+N^{\theta'}}\in B^{(2)}_x)\,\vert\, \mathcal{F}_{k_0} \Big] (\widetilde{\omega}) &> P^z_{\widetilde{\omega}}(X_M \in B_x^{(0)}) - CN^{-c\log N}
		\end{align}
		and $\bP(\widetilde{G}_3(z,N))\ge 1-CN^{-c\log N}$. The idea behind the last two inequalities is the following. Since $m_{k_0}=M+\log^2 N$ we can deduce that
		\begin{equation*}
			\bE\big[ P^z_\omega(X_M \in B_x^{(0)})\,\vert\, \cF_{k_0} \big](\widetilde{\omega}) \le P^z_{\widetilde{\omega}}(X_M\in B_x^{(0)})+ CN^{-c\log N},
		\end{equation*}
		holds with probability at least $1-C N^{-c \log N}$ by similar arguments as explained in the proof outline of Lemma~\ref{lem:BoundingU_k} below, see the discussion around the set $D(N)$ therein. The main idea is that the behaviour of the quenched random walk up to time $M$ mostly depends on the environment up to time at most $M+\log^2 (N)$, since the safety layer of time ensures that differences in $\xi$ up to time $M$ only happen with probability at most $CN^{-c\log N}$. Therefore it is sufficient to prove that the following claim holds: 
		\begin{align}
			\label{eq:21}
			\begin{split}
				\bE\Big[ P^z_\omega(X_{M+N^{\theta'}}\in B_x^{(1)}) &- P^z_\omega(X_M\in B_x^{(0)})\,\Big\vert\, \mathcal{F}_{k_0}\Big]\\ &\le \bE\Big[ P^z_\omega(X_{M+N^{\theta'}}\in B_x^{(1)},X_M\notin B_x^{(0)})\,\Big\vert\,\mathcal{F}_{k_0} \Big]\\
				&\le CN^{-c'\log N},
			\end{split}
		\end{align}
		for some $c'>0$. Indeed, combining this claim with the above inequalities (and analogous bounds in the other direction when we replace $B_x^{(1)}$ by $B_x^{(2)}$) we obtain that $\widetilde{G}_2(z,N)\cap\widetilde{G}_3(z,N)\subset G_1(z,\theta,N)$ and from which we conclude the statement of Lemma \ref{lem:Prop 3.11 Ersatz}.  		
		The remainder of this proof is therefore devoted to the proof of the above stated claim in \eqref{eq:21}. 
		For this,  for $(x',m')\in\bZ^d\times\bZ$ and $K\in \bN$,  
		let $A((x',m'),K)$ be the event that 
		\begin{equation*}
			P^{(x',m')}_\omega\big( \norm{X_n-x'} \ge \sqrt{n-m'}\log^3 (K)\big) \le CK^{-c\log(K)/2}
		\end{equation*}
		for each $n\in [m'+\tfrac{K}{2}, m'+K]$ and consider
		\begin{equation}
			A^{(x,M)} (N^{\theta'})= \bigcap_{x' \in B_{2N^{\theta'}}(x)} A((x',M),N^{\theta'}).
		\end{equation} 
		Then, using that $P^z_\omega(X_{M+N^{\theta'}}\in B_x^{(1)},X_M\notin B_x^{(0)}) \leq CN^{-c\theta'\log(N^{\theta'})/2}$ for any $\omega \in A^{(x,M)} (N^{\theta'})$, 
		we have that
		\begin{align*}
			\bE&\Big[ P^z_\omega(X_{M+N^{\theta'}}\in B_x^{(1)},X_M\notin B_x^{(0)})\,\Big\vert\, \mathcal{F}_{k_0}\Big]\\
			&\le\bE\Big[ \indset{ A^{(x,M)} (N^{\theta'}) }P^z_\omega(X_{M+N^{\theta'}}\in B_x^{(1)},X_M\notin B_x^{(0)})\,\Big\vert\,\mathcal{F}_{k_0} \Big] + \bE[\indset{ A^{(x,M)} (N^{\theta'})^\compl}\,\vert\, \mathcal{F}_{k_0}]
			\\&\le CN^{-c\theta'\log(N^{\theta'})/2} + \bE[\indset{ A^{(x,M)} (N^{\theta'})^\compl}\,\vert\, \mathcal{F}_{k_0}].
		\end{align*}
		To bound the  second term on the right hand side, we note first that, by the definition of $\mathcal{F}_{k_0}$ and the Markov property, 
		\begin{align*}
			&\bE[\indset{A^{(x,M)} (N^{\theta'})^\compl}\,\vert\, \mathcal{F}_{k_0}]\\
			&=\bP\Big( \exists\, n\in\big[M+\frac{N^{\theta'}}{2},M+N^{\theta'} \big]\colon P^{(x',M)}_\omega\big( \norm{X_n-x'} \ge \sqrt{n-M} \log^3(N^{\theta'}) \big) \ge CN^{-c\theta'\log(N^{\theta'})/2}\,\vert\, \mathcal{F}_{k_0} \Big) \\
			&\le \sum_{n=M+\frac{N^{\theta'}}{2}}^{M+N^{\theta'}} \sum_{y\in B_{\log^2 N}(x')} \bP\Big( P^{(y,M+\log^2 N)}_\omega\big( \norm{X_n-x'} \ge \sqrt{n-M} \log^3(N^{\theta'}) \big) \ge CN^{-c\theta'\log(N^{\theta'})/2}
			\Big).
		\end{align*}	
		Further, for any $y \in B_{\log^2 N}(x')$, by the Markov inequality it holds that
		\begin{align*}
			&\bP\Big( P^{(y,M+\log^2 N )}_\omega\big( \norm{X_n-x'} \ge \sqrt{n-M} \log^3(N^{\theta'}) \big) \ge CN^{-c\theta'\log(N^{\theta'})/2} \Big)\\
			&\le  \bP\Big( P^{(y,M+\log^2 N)}_\omega\big( \norm{X_n-y} \ge \sqrt{n-M} \log^3(N^{\theta'})-\log^2 N \big) \ge CN^{-c\theta'\log(N^{\theta'})/2} \Big)\\
			&\le  \frac{\bP^{(y,M+\log^2 N)}(\norm{X_n-y} \ge \sqrt{n-M} \log^3(N^{\theta'})-\log^2 N)}{CN^{-c\theta'\log(N^{\theta'})/2}}
		\end{align*}
		To obtain a suitable upper bound for the numerator of the last line we use the same arguments as in the proof of Lemma~\ref{lem:lem3.6SS} in Section~\ref{sec:pf_lem3.6SS}. The difference lies mainly in notation but for the readers convenience we provide the steps here also.
		Now, recall the event $R^\simu_N$ from \eqref{eq:defn_RN} and note that $\bP(R_N^\simu) \ge 1-CN^{-c\log N}$. 
		On the event $R^\simu_N$ the time between two consecutive regenerations is at most $\log^2 N$ up to time $N$. Moreover $(X_{T_k}-X_{T_{k-1}}, T_k-T_{k-1})_{1\le k\le N}$ conditioned on $R_N^\simu$ is an i.i.d.\ sequence and the increments $X_{T_k}-X_{T_{k-1}}$ are symmetrically distributed. Therefore, setting $\tilde{n}=n-M-\log^2(N)\ge \frac{N^{\theta'}}{2}-\log^2(N)$ and using translation invariance, we obtain
		\begin{align*}
			&\bP^{(y,M+\log^2 N)}(\norm{X_n-y} \ge \sqrt{n-M} \log^3(N^{\theta'})-\log^2 N)\\
			&=\bP^{(y,0)}(\norm{X_{\tilde{n}}-y} \ge \sqrt{\tilde{n}+\log^2 N} \log^3(N^{\theta'})-\log^2 N)\\
			&\le \bP^{(y,0)}(\norm{X_{\tilde{n}}-y} \ge \sqrt{\tilde{n}} \log^3(N^{\theta'})-\log^2 N \,\vert\, R^\simu_N) + CN^{-c\log N}\\
			&\le \bP^{(y,0)}(\exists k\le \tilde{n}\colon\norm{X_{T_k}-y} \ge \frac{1}{2}\sqrt{\tilde{n}} \log^3(N^{\theta'}) \,\vert\, R^\simu_N) + CN^{-c\log N}\\
			&\le \sum_{k=1}^{\tilde{n}} \bP^{(y,0)}(\norm{X_{T_k}-y}\ge \frac{1}{2}\sqrt{\tilde{n}}\log^3 (N^{\theta'})\,\vert\, R^\simu_N) + CN^{-c\log N}\\
			&\le d\sum_{k=1}^{\tilde{n}} \exp\bigg( -C\frac{\tilde{n}\log^6(N^{\theta'})}{k\log^4(N)} \bigg) +CN^{-c\log N}\\
			&\le CN^{-c\log N},
		\end{align*}
		where we used the fact that on $R^\simu_N$ the random walk cannot reach the distance $\sqrt{\tilde{n}}\log^3(N^{\theta'})-\log^2N$ and come back to $\tfrac{1}{2}\sqrt{\tilde{n}}\log^3(N)$ in between two regeneration times when $\tilde{n}\ge \frac{N^{\theta'}}{2}-\log^2(N)$ and Azuma's inequality for each coordinate in the fourth line. 	
		Since $\abs{B_{2N^{\theta'}}(x)}=(2N^{\theta'})^d$ and by translation invariance, we obtain the upper bound
		\begin{equation*}
			\bE[\indset{A^{(x,M)} (N^{\theta'})^\compl}\,\vert\, \mathcal{F}_{k_0}] \le 2^{d-1} N^{(d+1)\theta'}\log(N)^{2d} N^{-c(1-\frac{(\theta')^2}{2})\log(N)},
		\end{equation*}
		from which the claim is an immediate consequence.
	\end{proof}

	The proof of Lemma \ref{lem:BoundingU_k}  follows along the lines of Lemma~3.15 and Lemma~3.18 in \cite{SteibersPhD2017}. Although these lemmas  were stated in \cite{SteibersPhD2017} for $d\geq3$, their proofs  extend almost without changes to the cases when $d=1$ or $d=2$. For the brevity of this work we therefore present here only an outline of the argument to describe the ideas
	behind the proof and refer to \cite{SteibersPhD2017}  
	for the detailed proofs.
	\begin{proof}[Proof outline of Lemma \ref{lem:BoundingU_k}]
		Comparing $M_k$ and $M_{k-1}$, as defined in \eqref{eq:defn_martingale},  we see that the difference stems from knowing the value of $\omega$ on an additional site,  $z_k$. Therefore, to estimate $U_k$, we need to control what influence this knowledge has on the walk, i.e.\ we compare the two (random) probability laws 
		\begin{equation}
			\label{eq:twolaws}
			\bE \Big[ P^z_\omega (X \in \cdot )\,\vert\, \mathcal{F}_k \Big]
			\quad \text{and} \quad
			\bE \Big[ P^z_\omega (X \in \cdot )\,\vert\, \mathcal{F}_{k-1} \Big]
		\end{equation}
		on the set of random walk paths (recall that the starting point $z \in \widetilde{\cP}(N)$ is fixed in this subsection).
		
		For this we first note that, by Lemma~\ref{lem:lem3.6SS}, for any space-time starting point $z\in\widetilde{\cP}(N)$ the walk stays in $\cP(N)$ until time $N$ with probability at least $1-CN^{-c\log N}$. Thus it is sufficient to consider only the influence of $\omega(z_k)$ for $z_k \in \cP(N)$. We also restrict to $\omega$'s in the set $D(N)$ satisfying that, for some positive constant $\rho$ and each $(x,t) \in \cP(N)$:
		\begin{enumerate}
			\item either $(x,t) \notin \cC$ and the longest open path starting from it is at most $\log^2 N$ long;
			\item or $(x,t) \in \cC$ and it is connected via open paths to at least $\rho \log^2N$ many vertices at time $t+\log^2 N$.
		\end{enumerate}
		Note that, by classical estimates for supercritical oriented percolation, see e.g.\ Corollary~3.4 in \cite{SteibersPhD2017}, which hold in any dimension $d \ge 1$, we have $\bP(D(N))\ge 1-CN^{-c\log N}$.
		
		Under the above restrictions, the behaviour of the random walk under the two laws in \eqref{eq:twolaws} differs only if $\xi$ changes for some $(x, t) \in\cP(N)$ with $t \le m_k$ depending on the value of $\omega$ at $z_k$. To control this  we can consider two cases.
		
		Firstly, if the random walk $X$ does not visit the set $\Delta_k$, it is obvious, by the definition of $D(N)$, that a change of $\omega$ on $z_k$ leads to a change of $\xi$ in the areas that $X$ observes along its path with probability bounded by $CN^{-c\log N}$. Indeed, for any  $(x,t) \notin \Delta_k$ as above and $\omega \in D(N)$ we have the following dichotomy: if $\xi_t(x)=0$, then there is no open path connecting $(x,t)$ to $z_k$; if $\xi_t(x)=1$, then this value can only change if $\omega(z_k)=1$ is changed to $\omega(z_k)=0$ and every infinite open path starting from $(x,t)$ travels through $z_k$. But on $D(N)$, $(x,t) \in \mathcal{C}$ is connected to at least $\rho \log^2 N$ vertices at time $t + \log^2 N$. The probability that a contact process started from $\rho \log^2 N$ many vertices dies out is bounded by $CN^{-c\log N}$, see Lemma~1.3 in \cite{SteibersPhD2017} or Theorem~2.30 in \cite{Liggett:1999} for the continuous case. Similarly, the probability that the cluster starting from these $\rho \log^2 N$ many points at time $t + \log^2 N$ contains less than say $(\rho/2) \log^2 N$ points in time-slice $m_k$, is also bounded by $CN^{-c\log N}$.
		Since $\abs{\cP(N)}\le CN^{1+d/2}\log^{3d}N$ the probability that there exists a $(x,t) \notin \Delta_k$ for which $\xi$ changes its value if we flip the value of $\omega(z_k)$ is bounded by $CN^{-c\log N}$. Therefore, contributions from this case will be covered by the last term in \eqref{eq:BoundingU_k}.
		\smallskip
		
		The second case occurs when $X$ actually visits $\Delta_k$, which requires a more careful consideration and where we look at further two distinct possibilities, that $z\in\Delta_k$ or $z\notin \Delta_k$ particularly, as detailed next.
		
		For $z\in \Delta_k$, we note that there is a special regeneration point on the path of $X$ that is connected through open paths to all $z\in\Delta_k$ within at most $2\log^{6d+9}N$ steps with probability at least $1-CN^{-c\log N}$, see \cite[Lemma~3.16 and the discussions around it]{SteibersPhD2017} and note also the figure on p.~57 therein. (The proof given there works for any $d \ge 1$.) This allows us to instead compare the probabilities starting from that special regeneration point. Note that we then consider the annealed law since the special regeneration point lies beyond $z_k$ (in the time direction). Hence, up to an error term of $CN^{-c\log N}$, the estimate for $z\in\Delta_k$ comes down to a difference of annealed laws with starting positions in the set
		\begin{equation*}
			\boldsymbol{R}(z_k,N)\coloneqq\Big\{ (x,n)\colon \norm{y_k-x}\le 3\log^{6d+9} N, 0\le n-m_k\le 2\log^{6d+9} N \Big\}.
		\end{equation*}
		Thus for $z \in \Delta_k$ and $\omega \in D(N)$ we have the upper bound
		\begin{multline}
			\label{eq:13}
			U_k \le  CN^{-c\log N} + CP^z_\omega(\Delta_k\text{ is visited})\log^{6d+9}(N)\\
			\times \sup_{(x_i,n_i)\in\boldsymbol{R}(z_k,N)}\abs{\bP^{(x_1,n_1)}(X_{M+N^{\theta'}}=v)-\bP^{(x_2,n_2)}(X_{M+N^{\theta'}}=v)} 
		\end{multline}
		Therefore, using the annealed derivative estimates from Lemma~\ref{lem:annealed_derivative_estimates}, we obtain the upper bound $C\log^{6d+9}(N) N^{-\theta'\tfrac{d+1}{2}}$ in this case.
		
		For $z=(y, m) \notin \Delta_k$, note that $m_k-m>\log^2 N$, since we are in the case where $X$ visits $\Delta_k$. Furthermore we only have to consider the behaviour of the quenched law until $X$ hits $\Delta_k$ and can then, by the Markov property, refer to the previous part of this case, where we had $z\in \Delta_k$. Since the area in which the random walk moves until hitting $\Delta_k$ is separated by $\log^2 N$ from $z_k$ in the time coordinate, we again conclude that changes in the environment $\xi$ up to hitting $\Delta_k$ happen with probability at most $CN^{-c\log N}$. Thus, all in all, the upper bound for the case when $X$ visits $\Delta_k$ is given by
		$CP^z_\omega(\Delta_k\text{ is visited})\log^{6d+9}(N) \cdot N^{-\theta'\tfrac{d+1}{2}} + CN^{-c\log N}$. 
	\end{proof}

	\subsection{Proof of Theorem \ref{thm: Theorem 3.24 Ersatz}} \label{ssec:proof_Theorem_3.24_Ersatz}
	
	With Lemma \ref{lem:Prop 3.11 Ersatz} at hand, the proof of Theorem \ref{thm: Theorem 3.24 Ersatz} for $d\geq 1$, albeit somewhat long and technical, follows by the same reasoning as that for the proof of
	\cite[Theorem 3.24]{SteibersPhD2017}. We provide here a summary of the main steps, again building on  the so-called \emph{environment exposure} technique as in the proof of Lemma \ref{lem:Prop 3.11 Ersatz}.
	
	There are two additional ingredients we need for the proof of Theorem~\ref{thm: Theorem 3.24 Ersatz}. In essence these lead to an improvement of Lemma~\ref{lem:Prop 3.11 Ersatz}. The first one gives a bound on the quenched probability of hitting boxes. 
	
	\begin{lemma}[{Extension of \cite[Lemma 3.22]{SteibersPhD2017}}] 
		Let $d \ge 1$. For \label{lem:quenched_box_upper_bound} every $0<\theta\le 1$ and $z\in\widetilde{\cP}(N)$, let $G_2(z,\theta,h,N)$ be the event that for every $M 
		\in [\tfrac{2}{5}N,N]$ and every box $B$ of side length $N^{\theta/2}$ that is contained in $[-\sqrt{N}\log^3 (N),\sqrt{N}\log^3(N)]^d$, 
		\begin{equation*}
			P^z_\omega \big( X_M \in B \big) \le \log^h(N)N^{-d(1-\theta)/2}.
		\end{equation*}
		Then for every $0<\theta\le 1$ and $z\in\widetilde{\cP}(N)$ there exist constants $C,c>0$ and $h=h(\theta)\ge 0$, such that uniformly in $z$
		\begin{equation*}
			\bP\big(G_2(z,\theta,h,N)\big) \ge 1-CN^{-c\log N}.
		\end{equation*}
	\end{lemma}
	\begin{proof}
		A more detailed version of this proof can be found in \cite{SteibersPhD2017}, see Lemma~3.22 therein. We will restrict ourselves here to the core ideas. 
		The lemma is proved via induction over a decreasing sequence of $\theta$'s. Note that for $\theta>\tfrac{2d +1}{2d+2}$ the claim holds with $h=0$ by Lemma~\ref{lem:Prop 3.11 Ersatz}. Now assume the claim holds for some $0<\theta<1$ and fix a $\theta'$ with $\tfrac{2d+1}{2d+2}\theta<\theta'<\theta$. Furthermore,  set $\rho =\tfrac{\theta'}{\theta}>\frac{2d+1}{2d+2}$ and let $B =B_{N^{\theta'/2}}(x)$ with $x \in \bZ^d$. The set of all $\omega$ for which, writing $l(N)=2N^{\rho/2}\log^3(N)$, 
		\begin{multline*}
			\left| P^{z}_\omega\big( X_M\in B \big) -
			\sum_{v \in B_{l(N)}(x)} 
			P^{z}_\omega\big( X_{M-N^{\rho}}=v \big) P^{(v,M-N^\rho)}_\omega\big( X_M\in B \big) \right|
			\le CN^{-c\log N}
		\end{multline*}
		holds has probability at least $1-CN^{-c\log N}$ by applying Lemma~\ref{lem:lem3.6SS} to $P^{(v,M-N^\rho)}_\omega\big( X_M\in B_x \big)$ for $v$ not in the sum. Since $\rho>\tfrac{2d+1}{2d+2}$, for every box $B'$ with side length $N^{\rho/2}$ we have
		\begin{equation*}
			P^{z}_\omega\big( X_{M-N^\rho}\in B' \big)\le N^{-d(1-\rho)/2}
		\end{equation*}
		with probability at least $1-CN^{-c\log N}$, again by Lemma~\ref{lem:Prop 3.11 Ersatz}. Additionally, by induction assumption there exists $h'>0$ such that, with probability at least $1-CN^{-c\log N}$, 
		\begin{equation*}
			P^{(v,M-N^\rho)}\big( X_M\in B \big)\le \log^{h'}(N^\rho)(N^\rho)^{-d(1-\theta)/2}.
		\end{equation*}
		Note that here we are essentially considering events of the form $G_2(z,\theta,h',N^\rho)$, and thus boxes with side length $N^{\rho\theta/2}=N^{\theta'/2}$, i.e.\ the side length of $B$. 
		Lastly, the number of different boxes of side length $N^{\rho/2}$ needed to cover $x+[-N^{\rho/2}\log^3(N),N^{\rho/2}\log^3(N)]^d$ is bounded by $\log^{3d+1}(N)$.  
		Let $C(N)$ be one such cover. From this and the above bounds, we therefore have that 
		\begin{align*}
			&P^{z}_\omega\big( X_M\in B \big)\\ 
			&\le \sum_{B'\in C(N)}
			C P^{z}_\omega\big( X_{M-N^\rho}\in B' \big) \log^{h'}(N)N^{-d(\rho-\theta')/2} +CN^{-c\log N}\\
			&\le C\log^{3d+1}(N)N^{-d(1-\rho)/2}\log^{h'}(N)N^{-d(\rho-\theta')/2}\\
			&=C\log^{3d+h'+1}(N)N^{-d(1-\theta')/2},
		\end{align*}
		with probability at least $1-CN^{-c\log N}$, which concludes the proof since the bound holds irrespectively of the center $x$ of $B$. 
	\end{proof}
	A central part of the proof of Lemma~\ref{lem:Prop 3.11 Ersatz} is the upper bound on the expression
	\begin{equation*}
		\abs{\bE\big[ P^z_\omega(X_{N+N^{\theta'}}=v)\,\vert\, \cF_{k_0} \big] - \bP^z\big( X_{N+N^{\theta'}}=v \big)}
	\end{equation*}
	from Lemma~\ref{lem:McDiarmidApplied}, which was based on  a suitable upper bound on 
	$	\sum_{k=1}^{k_0} (P_\omega^z(\Delta_k \text{ visited}))^2$
	obtained in \eqref{eq:12}.
	The following lemma gives us an improved version of Lemma~\ref{lem:McDiarmidApplied} by estimating this sum more carefully, utilising the improved bounds obtained in Lemma~\ref{lem:quenched_box_upper_bound}.
	
	\begin{lemma}[{Replacement for \cite[Lemma 3.23]{SteibersPhD2017}}]
		\label{lem:better_comparison}
		Let $d\ge 1$ and recall $\cF_{k_0}$ from Lemma~\ref{lem:McDiarmidApplied}. Let $0<\eta<\tfrac{5}{4d}\wedge 1$ and define the event $G_3(z,N^\eta,N)$ as the event that for every $v\in\bZ^d$
		\begin{equation}
			\abs{\bE\big[ P^z_\omega(X_{N+N^{\eta}}=v)\,\vert\, \cF_{k_0} \big] - \bP^z\big( X_{N+N^{\eta}}=v \big)} \le N^{-d/2}N^{-\tfrac{d}{5}\eta}.
		\end{equation}
		Then there exist constants $C,c>0$ such that 
		\begin{equation*}
			\bP(G_3(z,N^\eta,N)) \ge 1-CN^{-c\log N}.
		\end{equation*}
	\end{lemma}
	\begin{proof}
		The main change in the approach compared with the proof of Lemma~\ref{lem:McDiarmidApplied} is that we divide $\cP(N)$ into ``slices'', which gives in the estimates more precise control on the number of remaining steps of the walk. 
		Particularly, let $L$ be a large integer such that $2^{-(L+1)}N \le N^\eta -\log^2 N < 2^{-L}N$ and, for $1\le l< L$, define
		\begin{equation*}
			\cP^{(l)} \coloneqq \cP(N) \cap\big\{ (x,n)\colon x\in\bZ^d, 2^{-l-1}N \le N-n < 2^{-l}N \big\}.
		\end{equation*}
		Additionally, set
		\begin{align*}
			\cP^{(0)}&\coloneqq \cP(N)\cap\big\{ (x,n)\colon x\in\bZ^d, 0\le n\le \tfrac{N}{2} \big\},\\
			\cP^{(L)}&\coloneqq \cP(N)\cap\big\{ (x,n)\colon x\in\bZ^d, 0\le N-n<2^{-L}N \big\},\\
			F(v)&\coloneqq \big\{ (x,n)\in\cP(N)\colon \norm{x-v}\le \sqrt{N+N^\eta-n} \cdot \log^3 N \big\},
		\end{align*}
		let $\cP^{(l)}(v)\coloneqq \cP^{(l)}\cap F(v)$, and put 
		\begin{equation*}
			V(l) \coloneqq \sum_{k\, \colon z_k\in\cP^{(l)}(v)}(P^z_\omega(\Delta_k \text{ visited}))^2.
		\end{equation*}
		
		As for Lemma~\ref{lem:McDiarmidApplied}, the proof is based on environment exposure and McDiarmid's inequality. 
		To obtain an improved bound on $U$  from \eqref{ineq:McDiarmid} compared to \eqref{eq:12}, we first derive bounds on the terms $V(l)$.
		
		We begin with the cases $d \in \{1, 2\}$. Recall the definition of $W(N)$ from \eqref{eq:defn_W(N)} and observe that, on $W(N)$, we have
		\begin{equation}
			\label{eq:V(0)_bound_d<=2}
			V(0) \le E^z_\omega\Big[ \#\big\{ n<N\colon \norm{X_n - X'_n} \le \log^2(N) \big\}\indset{R^\simu_N} \Big] \le N^{1/2 + \varepsilon}
		\end{equation}
		with a (small) $\varepsilon>0$ to be tuned later (recall that Lemma~\ref{lem:replacement_cor_3.20} works for every $\varepsilon>0$).        
		Further, recall the space-time set $\Delta_k$ from \eqref{eq:Delta_k} and denote its bottom boundary by
		\begin{equation*}
			\partial^-\Delta_k = \{ (x,n) \in \Delta_k\colon n=\min\{ i \colon (x,i) \in \Delta_k \} \}.
		\end{equation*}
		For $l\geq 1$, by Lemma~\ref{lem:quenched_box_upper_bound}, we have that, with probability at least $1-CN^{-c\log N}$,
		\begin{align*}
			V(l) &= \sum_{k\colon z_k\in\cP^{(l)}(v)} \Big( P^z_\omega(X_{m_k-\log^2(N)}\in \partial^-\Delta_k) \Big)^2\\
			&\le \abs{\{k\colon z_k\in\cP^{(l)}(v)\}} \big(\log^h(N)N^{-d(1-\theta)/2}\big)^2\\
			&\le C2^{-l\tfrac{2+d}{2}}N^{1+\tfrac{d}{2}}\log^{2h+3d}(N)N^{-d(1-\theta)}
		\end{align*}
		for a (small) $\theta>0$ to be tuned later. Note that $h$ (from Lemma~\ref{lem:quenched_box_upper_bound}) depends on the choice of $\theta$ and we used that the side length of $\partial^-\Delta_k$ is $2\log^2 N$ (which is less than $N^\theta$ for any $\theta > 0$ when $N$ is large enough).
		Recall from \eqref{eq:13} the upper bound
		\begin{multline*}
			U_k \le  CN^{-c\log N} + CP^z_\omega(\Delta_k\text{ is visited})\log^{6d+9}(N)\\
			\times \sup_{(x_i,n_i)\in\boldsymbol{R}(z_k,N)}\abs{\bP^{(x_1,n_1)}(X_{N+N^{\eta}}=v)-\bP^{(x_2,n_2)}(X_{N+N^{\eta}}=v)}.
		\end{multline*}
		For $z_k\in\cP^{(l)}(v)$ we then have
		\begin{equation*}
			\sup_{(x_i,n_i)\in\boldsymbol{R}(z_k,N)}\abs{\bP^{(x_1,n_1)}(X_{N+N^{\eta}}=v)-\bP^{(x_2,n_2)}(X_{N+N^{\eta}}=v)} \le C\log^{6d+9}(N)(2^{-l}N)^{-\tfrac{d+1}{2}}.
		\end{equation*}
		Combining the estimates obtained above we conclude 
		that on the event $W(N)$
		\begin{align}
			\frac{U}{\log^{24d+36}(N)}
			&\le \sum_{l=0}^L C(2^{-l}N)^{-(d+1)}V(l) + CN^{-c\log N} \notag \\
			&\le N^{-(d+1)}V(0) + \sum_{l=1}^L C(2^{-l}N)^{-(d+1)}2^{-l\tfrac{2+d}{2}}N^{1+\tfrac{d}{2}}\log^{2h+3d}(N)N^{-d(1-\theta)}\notag \\
			&\le N^{\varepsilon-d-1/2} + C\log^{2h+3d}(N)N^{-d\tfrac{3}{2}+d\theta}\sum_{l=1}^L 2^{l(d+1)}2^{-l\tfrac{2+d}{2}}\notag \\
			&\le N^{\varepsilon-d-1/2} + C\log^{2h+3d}(N)N^{-d\tfrac{3}{2}+d\theta}2^{L\tfrac{d}{2}} \notag \\
			&\le N^{\varepsilon-d-1/2} + C\log^{2h+3d}(N)N^{-d\tfrac{3}{2}+d\theta}N^{(1-\eta)d/2} \notag \\
			&\le N^{\varepsilon-d-1/2} + C\log^{2h+3d}(N)N^{-d(1-\theta)}N^{-\eta\tfrac{d}{2}}\notag \\
			&\le C\log^{2h+3d}(N)N^{-d(1-\theta)}N^{-\eta\tfrac{d}{2}},
			\label{eq:pf.better_comparison.2}
		\end{align}
		for
		\begin{equation}
			\label{eq:pf.better_comparison.1}
			0 < \varepsilon\le \tfrac{1}{2}+d\theta-\eta\tfrac{d}{2}.
		\end{equation}
		Now the assertion follows using McDiarmid's inequality, for $\theta < \tfrac{1}{10}\eta$.
		Note that choices for $\varepsilon>0$ and $0 < \theta < \tfrac{1}{10}\eta$ which are compatible with
		\eqref{eq:pf.better_comparison.1} are possible whenever $0<\eta<\tfrac{5}{4d}$ (and the argument breaks down for $\eta \ge 1$), which is possible when $d=1$ or $2$.
		
		The case $d \ge 3$ is analogous (and also covered by \cite[Lemma~3.23]{SteibersPhD2017}) except that
		we can replace \eqref{eq:V(0)_bound_d<=2} by the much tighter bound 
		\begin{equation}
			\label{eq:V0_bound_dge3}
			V(0)\le \log^{12}(N),
		\end{equation}
		see \cite[Corollary~3.20]{SteibersPhD2017},     
		owing to the fact that random walkers in $d\ge 3$ are transient in our model.
		This means that in the derivation of the bound on $U$ leading to \eqref{eq:pf.better_comparison.2}
		we can replace $\varepsilon$ by $-1/2$, leading to
		\begin{align*}
			U &\le  \log^{12}(N) N^{-(d+1)} + C\log^h(N)N^{-d(1-\theta)}N^{-\eta\tfrac{d}{2}}\\
			&\le C(\log N)^{\max(h,12)}N^{-d(1-\theta)}N^{-\eta\tfrac{d}{2}},
		\end{align*}
		if $0<1+d\theta-\eta\tfrac{d}{2}$. A compatible choice for $0 < \theta < \tfrac{1}{10}\eta$
		requires in fact (only) $0<\eta<\tfrac{5}{2d}$ in this case.
	\end{proof}

	\begin{proof}[Proof of Theorem \ref{thm: Theorem 3.24 Ersatz}] With the changes discussed above, the proof follows the template from \cite[Theorem~3.24]{SteibersPhD2017}.
		
		Fix $0<\theta\le 1$ and $z\in\widetilde{\cP}(N)$. Let $\tfrac{3}{16}\theta<\theta'<\tfrac{7}{20}\theta$ and set $V\coloneqq N^{\tfrac{2\theta'}{d}}$. By Lemma~\ref{lem:better_comparison} we know that there exists a set $G_3=G_3(z,V,N)$ such that $\bP(G_3)\ge 1-CN^{-c\log N}$ and that on $G_3$ it holds that, for any $v\in\bZ^d$,
		\begin{equation}
			\label{eq:20}
			\abs{\bE\big[ P^z_\omega(X_{N+V}=v)\,\vert\, \cF_{k_0} \big] - \bP^z(X_{N+V}=v)} \le N^{-d/2}V^{-d/5}.
		\end{equation} 
		
		For $x\in[-N \log^3(N),N \log^3(N)]^d$ write $B_x^{(0)}=B_{N^{\theta/2}}(x)$. 
		Moreover, similar to the proof of Lemma~\ref{lem:Prop 3.11 Ersatz}, we consider boxes 
		\[ B_x^{(1)} = B_{N^{\theta/2}-\log^3(N)\sqrt{V}}(x) \quad \text{ and } \quad  B_x^{(2)}=B_{N^{\theta/2}+\log^3(N)\sqrt{V}}(x).\] 
		Note that, since $\tfrac{2\theta'}{d} < \theta$, we have 
		\begin{align*}
			\abs{B_x^{(2)}\setminus B_x^{(1)}} &\le \big( N^{\theta/2}+\log^3(N)\sqrt{V} \big)^d - \big( N^{\theta/2}-\log^3(N)\sqrt{V} \big)^d\\
			&\le CN^{(d-1)\theta/2}\log^{3}(N)\sqrt{V}.
		\end{align*} 
		and, moreover, 
		\begin{align*}
			\abs{B_x^{(2)}} &\le (N^{\theta/2}+\log^3(N)\sqrt{V})^d\\
			&=(N^{\theta/2}+\log^3(N)N^{\tfrac{\theta'}{d}})^d
			\le CN^{d\theta/2}.
		\end{align*}
		Combining the above, we have that there exists a set $G_4\subset G_3$ such that for all $\omega \in G_4$
		\begin{align*}
			\abs{ P^z_\omega&(X_N\in B_x^{(0)}) - \bP^z(X_N\in B_x^{(0)}) }\\
			&\le \abs{ \bE\big[ P^z_{\widetilde{\omega}}(X_{N+V}\in B_x^{(2)})\,\vert\,\cF_{k_0} \big] - \bP^z(X_{N+V}\in B_x^{(1)}) } + CN^{-c\log N}\\
			&\le \abs{\bE\big[ P^z_{\widetilde{\omega}}(X_{N+V}\in B_x^{(2)})\,\vert\,\cF_{k_0} \big] - \bP^z(X_{N+V}\in B_x^{(2)})} + \bP^z(X_{N+V}\in B_x^{(2)}\setminus B_x^{(1)})\\
			&\hspace{0.5cm} + CN^{-c\log N}\\
			&\le \abs{B_x^{(2)}} N^{-d/2}V^{-d/5} + \abs{B_x^{(2)}\setminus B_x^{(1)}} N^{-d/2} + CN^{-c\log N}\\
			&\le CN^{-d(1-\theta)/2}V^{-d/5} +CN^{(d-1)\theta/2}\log^{3}(N)V^{1/2} N^{-d/2} + CN^{-c\log N}\\
			&\le CN^{-\tfrac{d(1-\theta)}{2} - \tfrac{2}{5}\theta'}  + C\log^3(N)N^{-\tfrac{d(1-\theta)}{2} - \tfrac{\theta}{2}+\tfrac{\theta'}{d}} + CN^{-c\log(N)}\\
			&\le N^{-d(1-\theta)/2} N^{-\tfrac{3}{20}\theta},
		\end{align*}
		and $\bP(G_4) \ge 1-CN^{-c\log N}$. 
		Here we used an adaptation of equations~\eqref{eq:14} through \eqref{eq:18} in the first inequality. Note that this adaptation can be proved in the same way using Lemma~\ref{lem:lem3.6SS}. The third inequality uses \eqref{eq:20} and \eqref{eq:ann_der_est_proof_1}. The rest then follows by the estimates obtained above.
		For the last line we want to highlight that by the specific choice of $\tfrac{3}{16}\theta<\theta'<\tfrac{7}{20}\theta$ we have $-\tfrac{2}{5}\theta' < -\tfrac{3}{40}\theta$ and $-\tfrac{\theta}{2}+\tfrac{\theta'}{d} < -\theta\tfrac{10d-7}{20d}$. Therefore we can drop the constant $C$ from the previous lines.
	\end{proof}

	\section{Proofs of the main results}\label{sec:proofs}
	
	In this section we outline the proof of Theorem \ref{thm:main1} and
	Theorem \ref{thm:main2}, following the general exposition provided in
	\cite[Section 2]{BethuelsenBirknerDepperschmidtSchluter} closely,
	which in turn was highly influenced by the methodology of
	\cite{BergerCohenRosenthal2016}. Provided with the results from the
	previous section, these proofs for $d=1$ and $d=2$ require almost no
	changes compared to the proofs of their equivalent statements in
	\cite[Section 2]{BethuelsenBirknerDepperschmidtSchluter} for $d\geq3$
	and we therefore only highlight  the modifications needed.

	\subsection{An improved quenched annealed-comparison result}
	
	As a first consequence of Theorem~\ref{thm: Theorem 3.24 Ersatz} we
	conclude a further result, extending
	\cite[Lemma~2.1]{BethuelsenBirknerDepperschmidtSchluter} to $d=1,2$ by
	a slight modifications of the proof. To state this precisely, for $N,M
	\in \bN$ and $c, C>0$, let $K(N) \coloneqq K(N,M,c,C)$ denote the set
	of environments $\omega \in \Omega$ such that, for every $x\in B_N$, 
	\begin{align}
		\label{eq:claim1}
		\sum_{B \in \Pi_{M}} \bigl|P_{\omega}^{(x,0)} (X_{N} \in B ) -
		\bP^{(x,0)}(X_{N} \in B) \bigr| \le \frac{C}{M^{c}} +
		\frac{C}{N^{c}}.
	\end{align}
	Here, and later, for any positive real number $L$, we denote by
	$\Pi_L$ a partition of $\bZ^d$ into boxes of side length $\lfloor L
	\rfloor$.
	
	\begin{lemma}[Extension of Lemma 2.1 in \cite{BethuelsenBirknerDepperschmidtSchluter}]
		Let \label{claim:1 for all dimensions} $d \ge 1$ and $M \in \bN$. If
		$c>0$ is small enough and $C<\infty$ large enough, there are
		universal positive constants $\tilde c$, $\widetilde C$, for which
		we have
		\begin{align}
			\label{eq:78}
			\bP\bigl(K(N)\bigr)\ge 1-\widetilde CN^{-\tilde c\log N}
			\quad \text{for all } N.
		\end{align}
	\end{lemma}

	\begin{proof}[Proof outline]
		As mentioned in \cite[Section
		4]{BethuelsenBirknerDepperschmidtSchluter}, the proof is a
		``consequence of Proposition 4.1 (and Theorem 8.1) and can be proven
		analogously to the argument in the last part of the proof of Theorem
		5.1 in [2], page 2920." This part of the argument does not rely on
		any dimension dependent quantities and goes through verbatim. The
		rather lengthy proof of \cite[Proposition
		4.1]{BethuelsenBirknerDepperschmidtSchluter}, provided therein in
		Section~8, and partly in Appendix~C, hinges on Theorem~8.1 there, which we
		concluded extends to $d\geq 1$ in Theorem \ref{thm: Theorem 3.24
			Ersatz}. In addition, the arguments make use of the annealed
		derivative estimates of \cite[Lemmas 3.1 and
		3.2]{BethuelsenBirknerDepperschmidtSchluter} that extend to $d=1$
		and $d=2$ as concluded in Lemma
		\ref{lem:annealed_derivative_estimates} and Lemma
		\ref{lem:additional_annealed_estimate}, respectively, and 
		estimates for the discrete-time contact
		process e.g.\ from \cite[Lemma 8.5]{BethuelsenBirknerDepperschmidtSchluter}
		which hold in any dimension $d \ge 1$.
	\end{proof}

	\subsection{Existence of a measure absolutely continuous with respect to the environment: proof of Theorem \ref{thm:main1}}

	We are now in position to  argue that the following extension of  \cite[Proposition 2.2 and Corollary 2.4]{BethuelsenBirknerDepperschmidtSchluter} holds. To that end, recall the shift-operator $\sigma$ defined below \eqref{eq:15}. For $N\in \bN$, we define $Q_N$ by
	\begin{align}
		\label{eq:defn_Q_N}
		Q_N(A) \coloneqq  \bE\Big[\sum_{x \in \bZ^d}
		P^{(o,0)}_\omega(X_N =x)\ind{\sigma_{(x,N)} \omega \in A} \Big].
	\end{align}
	and consider the Ces\`aro
	sequence
	\begin{align}
		\label{eq:80}
		\widetilde{Q}_n \coloneqq \frac{1}{n}\sum_{N=0}^{n-1}Q_N, \quad n =1,2,\dots.
	\end{align}
	
	Since $( \widetilde{Q}_n )$ are measures on a compact space
	(using that $\Omega$ carries the product topology), the
	sequence is tight. In particular, there is a weakly converging
	subsequence, say $(\widetilde{Q}_{n_k})_k$, and we define its
	weak limit by 
	\begin{align}
		\label{eq:77}
		Q \coloneqq  \lim_{k\to \infty}\widetilde{Q}_{n_k}.
	\end{align}

	\begin{lemma}[Extension of Proposition 2.2 and Corollary 2.4 in \cite{BethuelsenBirknerDepperschmidtSchluter}]
		Let $d\ge 1$. \label{lem:con_prop_Q for all dimensions} There exists a
		universal constant $c>0$ so that for every $\varepsilon >0$ there is
		$M_0=M_0(\varepsilon) \in \bN$ and for every $M \ge M_0$ we have
		\begin{align}
			&\label{eq:29}
			\bP\Bigl(\Big\lvert\frac{1}{\abs{B_M}}
			\sum_{x\in B_M}\frac{dQ_N}{d\bP}(\sigma_{(x,0)}\omega)-1
			\Big\rvert>\varepsilon \Bigr) \le M^{-c\log M}.
			\\       &\label{eq:29a}
			\bP\Bigl(\Big\lvert\frac{1}{\abs{B_M}}
			\sum_{x\in B_M}\frac{dQ}{d\bP}(\sigma_{(x,0)}\omega)-1
			\Big\rvert>\varepsilon \Bigr) \le M^{-c\log M}.
		\end{align}
	\end{lemma}
	
	\begin{proof}[Proof outline]
		The arguments of \cite[Section
		5]{BethuelsenBirknerDepperschmidtSchluter}, albeit lengthy,
		go through word-for-word up to the proof of Corollary 2.4
		therein, only replacing the use of Lemma 2.1 with Lemma
		\ref{claim:1 for all dimensions}, and yields the claimed
		inequalities.
	\end{proof}
	
	Finally, based on the above, we have all the tools available
	in order to conclude Theorem \ref{thm:main1}.
	
	\begin{proof}[Proof of Theorem \ref{thm:main1}]
		By construction and shift invariance of $\bP$ we have $Q_N\ll \bP$
		for every $N$ and therefore $\tilde Q_n \ll \bP$ for every $n$.
		Furthermore, by exactly the same arguments as on Page 8 in
		\cite{BethuelsenBirknerDepperschmidtSchluter}, Equation 2.8 therein
		follows as a consequence of inequality \eqref{eq:29}, that is,
		\begin{align}
			\label{eq:82}
			\begin{split}
				\bE\Bigl[ \Bigl( \frac{1}{n}\sum_{N=0}^{n-1}\frac{dQ_N}{d\bP}
				\Bigr)^2 \Bigr]
				& = \frac{1}{n^2} \sum_{N,N'=0}^{n-1} \bE\Big[
				\frac{dQ_N}{d\bP}\frac{dQ_{N'}}{d\bP} \Big]
				\le
				\tilde c.
			\end{split}
		\end{align}
		From this we conclude that the family of Radon--Nikodym
		derivatives $(d\tilde Q_n/d\bP)_{n=1,2,\dots}$ is uniformly
		integrable. These facts together imply that we also have $Q \ll \bP$
		for any $Q$ obtained as in \eqref{eq:77}. A standard argument shows that $Q$ is invariant with respect to the
		point of view of the particle; see Proposition~1.8 in
		\cite{LiggettIPS1985} for an abstract argument or the proof of Lemma~1
		in \cite{DrewitzRamirez2014} for the argument in the case of random
		walks in random environments. 
		Lastly, for uniqueness of $Q$,   see \cite[Remark~2.6]{BethuelsenBirknerDepperschmidtSchluter}, whose validity transfers to $d=1$ and $d=2$ without any further ado.
	\end{proof}

	\subsection{The local limit theorem: proof of Theorem \ref{thm:main2}} 
	
	The proof of Theorem \ref{thm:main2} hinges on  the extension of \cite[Proposition 2.9]{BethuelsenBirknerDepperschmidtSchluter}  to $d=1$ and $d=2$ that we present in the following. For this, we first recall the definition of the hybrid measures form \cite{BethuelsenBirknerDepperschmidtSchluter}. 
	
	\begin{definition}[Definition 2.7 in  \cite{BethuelsenBirknerDepperschmidtSchluter}]
		For $\omega \in \Omega$ \label{defn:auxiliary_measures} and a given partition $\Pi(\ell)$ of $\bZ^d$ into
		boxes of a fixed side length $\ell \in \bN$, let 
		\begin{align}
			\label{eq:measure_ann_by_prefactor}
			\nu^{\annpre}_\omega(x,n)
			& 
			\coloneqq \frac{1}{Z_{\omega,n}}\bP^{(o,0)}(X_n=x)\varphi(\sigma_{(x,n)}\omega),\\
			\label{eq:measure_quenched}
			\nu^{\que}_\omega(x,n)
			& 
			\coloneqq  P_\omega^{(o,0)}(X_n=x),\\
			\label{eq:measure_box-quenched_by_pre}
			\nu^{\boxpre}_{\omega}(x,n)
			&
			\coloneqq P_\omega^{(o,0)}(X_n \in B^{(x)})
			\frac{\varphi(\sigma_{(x,n)}\omega)}{\sum_{y\in B^{(x)}}\varphi(\sigma_{(y,n)}\omega)}.
		\end{align}
		These are measures on
		$\bZ^{d+1}$, where 
		$Z_{\omega,n}=\sum_{x\in\bZ^d}\bP^{(o,0)}(X_n=x)\varphi(\sigma_{(x,n)}\omega)$
		is the normalizing constant in \eqref{eq:measure_ann_by_prefactor}
		and $B^{(x)}$ in \eqref{eq:measure_box-quenched_by_pre} is the
		unique $d$-dimensional box that contains $x$ in the partition $\Pi(\ell)$.
	\end{definition}
	
	For the normalization factor we have the following.
	\begin{proposition}[Extension of Proposition~2.8 from \cite{BethuelsenBirknerDepperschmidtSchluter}]
		\label{lem:limit_Z_omega for all dimensions}
		Let $d\ge 1$. For $\bP$-almost all $\omega\in\Omega$ the
		normalizing constant $Z_{\omega,n}$ satisfies
		\begin{align}
			\label{eq:81}
			\lim_{n\to\infty} Z_{\omega,n} =1.
		\end{align}
	\end{proposition}
	\begin{proof}[Proof outline]
		Proposition \ref{lem:limit_Z_omega for all dimensions} can be proven analogously to the proof of Proposition~2.8 in \cite{BethuelsenBirknerDepperschmidtSchluter}, see Section 6 therein. For this, one needs to replace the use of \cite[Corollary 2.4]{BethuelsenBirknerDepperschmidtSchluter} by Lemma \ref{lem:con_prop_Q for all dimensions}. Moreover, since the annealed CLT from \cite{BirknerCernyDepperschmidtGantert2013} (referred to directly after \cite[Equation 6.1]{BethuelsenBirknerDepperschmidtSchluter}) holds for $d\geq 1$, also the corresponding estimate is valid. Similarly, all annealed estimates based on  \cite[Lemma~3.1 and  Lemma~3.2]{BethuelsenBirknerDepperschmidtSchluter} should be replaced by Lemma \ref{lem:annealed_derivative_estimates} and Lemma \ref{lem:additional_annealed_estimate}, respectively. The remaining parts of the proof in \cite[Section 6]{BethuelsenBirknerDepperschmidtSchluter} should be kept unchanged, thereby yielding its extension to $d=1$ and $d=2$.
	\end{proof}
	
	The following proposition is the key ingredient for the proof of
	Theorem~\ref{thm:main2} and states that for large $n$ the above
	introduced measures are close to each other in a suitable norm. To
	state this precisely, for $\omega \in \Omega$ and any two probability
	measures $\nu^1_\omega$ and $\nu^2_\omega$ on $\bZ^d \times \bZ$ (more
	precisely these are transition kernels from $\Omega$ to
	$\bZ^d \times \bZ$) let the \emph{$L^1$ distance of $\nu^1_\omega$ and
		$\nu^2_\omega$} at time $n\in\bZ$ be defined by
	\begin{align}
		\label{eq:40}
		\norm{\nu^1_\omega-\nu^2_\omega}_{1,n} \coloneqq \sum_{x\in\bZ^d}
		\abs{\nu^1_\omega(x,n)-\nu^2_\omega(x,n)}.
	\end{align}
	Furthermore, for $k\le n$ the \emph{space-time convolution of
		$\nu^1_\omega$ and $\nu^2_\omega$} is defined by
	\begin{align}
		\label{eq:41}
		(\nu^1\ast \nu^2)_{\omega,k}(x,n) \coloneqq
		\sum_{y\in\bZ^{d}}\nu^1_\omega(y,n-k)\nu^2_{\sigma_{(y,n-k)}\omega}(x-y,k).
	\end{align}

	\begin{proposition}[Extension of Proposition 2.9 in \cite{BethuelsenBirknerDepperschmidtSchluter}]
		Fix \label{lem:main for all dimensions} $0<2\delta<\varepsilon<\frac{1}{4}$, and for
		$n\in\bN$ set $k=\lceil n^\varepsilon\rceil$ and
		$\ell=\lceil n^\delta \rceil$.  With respect to the partition  $\Pi(\ell)$, for $\bP$-almost every
		$\omega\in \Omega$, the measures from
		Definition~\ref{defn:auxiliary_measures} satisfy
		\begin{align}
			\label{eq:ann_to_ann-quenched}
			\tag{L1}
			\lim_{n\to\infty}
			\norm{\nu^{\ann\times\pre}_\omega-(\nu^{\ann\times\pre}\ast
				\nu^{\que})_{\omega,k}}_{1,n} & =0,\\
			\label{eq:ann-quenched_to_box-quenched}
			\tag{L2}
			\lim_{n\to \infty}\norm{(\nu^{\ann\times\pre}\ast
				\nu^{\que})_{\omega,k}-(\nu^{\mathrm{box-que}\times\pre}\ast\nu^{\que})_{\omega,k}}_{1,n}
			&  =0,\\
			\label{eq:box-quenched_to_quenched}
			\tag{L3}
			\lim_{n \to \infty}
			\norm{(\nu^{\mathrm{box-que}\times\pre}\ast\nu^{\que})_{\omega,k}
				-(\nu^{\que}\ast\nu^{\que})_{\omega,k}}_{1,n} & =0.
		\end{align}
	\end{proposition}

	\begin{proof}[Proof outline]
		This can be proven analogously to the proof of Proposition 2.9 in \cite{BethuelsenBirknerDepperschmidtSchluter}, see Section 7 therein. The proof is separated into three parts, each part concerning one of the inequalities \eqref{eq:ann_to_ann-quenched}--\eqref{eq:box-quenched_to_quenched}. We argue now that the techniques used to prove those three inequalities are independent of the dimension:
		\begin{itemize}
			\item[\eqref{eq:ann_to_ann-quenched}] The proof provided in \cite[Section 7]{BethuelsenBirknerDepperschmidtSchluter} is based on Lemma~3.1, Proposition~5.4, Corollary~2.4 and Proposition~2.8 from \cite{BethuelsenBirknerDepperschmidtSchluter}, and Lemma~3.6 from \cite{SteibersPhD2017}. The first and last of these holds for all $d\geq 1$ as noted in Lemma \ref{lem:annealed_derivative_estimates} and Lemma \ref{lem:additional_annealed_estimate}, respectively. Moreover, we concluded  in Lemma~\ref{lem:con_prop_Q for all dimensions} and Proposition~\ref{lem:limit_Z_omega for all dimensions} that Corollary~2.4 and Proposition~2.8 of \cite{BethuelsenBirknerDepperschmidtSchluter} also hold for all $d\ge 1$. Lastly, Proposition~5.4 \cite{BethuelsenBirknerDepperschmidtSchluter} states that, for $\bP$-almost every $\omega$, every $n \in \bN_0$, every $x \in \bZ^d$ and all $k \leq n$, it holds that $\phi(\sigma_{(x,n)} \omega) = \sum_{y \in \bZ^d} P_{\omega}^{(x+y,n-k)} (X_n =x) \phi(\sigma_{(x+y,n-k)}\omega),$ 
			whose proof follows by a general and dimension-independent argument, analogous to that of  \cite[Proposition 7.1]{BergerCohenRosenthal2016}.       
			\item[\eqref{eq:ann-quenched_to_box-quenched}] The proof provided in \cite[Section 7]{BethuelsenBirknerDepperschmidtSchluter} is based on Lemma~2.1, Corollary~2.4, Proposition~2.8, Lemma~3.1 and Theorem~8.1 from \cite{BethuelsenBirknerDepperschmidtSchluter}. The respective results for $d\ge 1$ are, in order, Lemma~\ref{claim:1 for all dimensions}, Lemma~\ref{lem:con_prop_Q for all dimensions}, Lemma \ref{lem:annealed_derivative_estimates}, Lemma~\ref{lem:limit_Z_omega for all dimensions} and Theorem~\ref{thm: Theorem 3.24 Ersatz}.
			\item[\eqref{eq:box-quenched_to_quenched}] The main tool to prove this inequality for $d\ge 3$ is Lemma~7.1 from \cite{BethuelsenBirknerDepperschmidtSchluter}, proven in Section 9 therein. The proof of this lemma relies on Lemma 2.1, Lemma A.1, Lemma 3.6, and Theorem 1.1 of \cite{BethuelsenBirknerDepperschmidtSchluter} and Lemma 3.6 from \cite{SteibersPhD2017}. Here, \cite[Theorem 1.1]{BethuelsenBirknerDepperschmidtSchluter} is the annealed local central limit theorem which was proven therein to hold for $d\geq 1$. Moreover, \cite[Lemma A.1]{BethuelsenBirknerDepperschmidtSchluter}  provides quantitative bounds for the discrete-time contact process, which holds true for $d\geq 1$ (literately,  it is stated in \cite[Lemma A.1]{BethuelsenBirknerDepperschmidtSchluter} to hold for $d\geq 2$ following the convention of \cite{GrimmettHiermerDPRW2000}, which corresponds to $d\geq 1$ in the notation of this paper). Therefore, since all the technical lemmas needed apply for $d\geq1$, Lemma 7.1 and so also the claim \eqref{eq:box-quenched_to_quenched}, holds in all dimensions.  
		\end{itemize}
	\end{proof}
	
	\begin{proof}[Proof of  Theorem \ref{thm:main2}]
		With Propositions \ref{lem:limit_Z_omega for all dimensions} and  \ref{lem:main for all dimensions} at hand for $d\geq1$, we conclude the proof ad verbum as in the proof of Theorem 1.4 in \cite{BethuelsenBirknerDepperschmidtSchluter}.
	\end{proof}


	\begin{appendix}
		\section{Proofs of preliminary annealed estimates}\label{sec:ppae}

		\subsection{Proof of Lemma \ref{lem:annealed_derivative_estimates}}
		\label{sec:pf_lem:annealed_derivative_estimates}
		A version of Lemma~\ref{lem:annealed_derivative_estimates} for $d\ge 3$ can be found as Lemma~3.9 in \cite{SteibersPhD2017} and was proved in the appendix therein. In this section we argue  along the same lines that these results also extend to $d=1,2$. We write an abbreviated version of the proof and refer for a more detailed version to \cite{SteibersPhD2017}, but nevertheless we hope to convey the core idea here. First we recall two lemmas that are used in the proof. 
		
		\begin{lemma}[Lemma~A.2 in \cite{SteibersPhD2017}]
			Let $d\ge 1$, let $\{Y_i\}_{i\ge 1}$ be a sequence of $\bZ^d$-valued random variables and $\{Z_i\}_{i\ge 1}$ a sequence of $\bN$-valued random variables, such that $\{(Y_i,Z_i)\}_{i\ge 1}$ are i.i.d.\ with respect to some probability measure $P$ with $\bE[\norm{Y_1}^3+Z_1^3]<\infty$. Assume in addition that there exists $v\in\bZ^d$, $k\in\bN$ such that $P((Y_1,Z_1)=(v,k))>0$ and $P((Y_1,Z_1)=(w,k+1))>0$ for every $w$ with $\norm{w-v}\le 1$. Let $S_n=\sum_{i=1}^n Y_i$ and $T_n=\sum_{i=1}^n Z_i$. Then there exists $C<\infty$ 
			such that for every $n,m\in\bN$ and every $x,y\in\bZ^d$ with $\norm{x-y}=1$
			\begin{align}
				\label{eq:lemmaA2_1}
				P((S_n,T_n)=(x,m)) \le Cn^{-\frac{d+1}{2}},
			\end{align}
			\begin{align}
				\label{eq:lemmaA2_2}
				\abs{P((S_n,T_n)=(x,m))-P((S_n,T_n)=(x,m+1))} \le Cn^{-\frac{d+2}{2}},
			\end{align}
			and
			\begin{equation}
				\label{eq:lemmaA2_3}
				\abs{P((S_n,T_n)=(x,m))-P((S_n,T_n)=(y,m))} \le Cn^{-\frac{d+2}{2}}.
			\end{equation}
		\end{lemma}
		\begin{proof} This follows by local CLT estimates obtained for example in \cite[Theorems\ 2.4.5 and 2.3.6]{LawlerLimic2010}. 
		\end{proof}
		
		For the next lemma, recall the sequence of regeneration times $(T_i)_i$ for a single random walker obtained in \cite{BirknerCernyDepperschmidtGantert2013}, see Section~2.2 and more specifically equation~(2.9) therein. We refer also to the discussion in Remark~\ref{rem:ontheclusterornot?} of the present text and the references given there.
		Note that in our setting, under $\bP^{(y,m)}$  the random walk starts at $(y,m)$ and we therefore let the ``first'' regeneration time $T_0$ be a random variable with $\bP^{(y,m)}(T_0=m)=1$.
		Lastly, for $\ell \leq M$ in $\bN$ and $m\le k \le M-\ell-1$, consider the events
		\begin{align}
			\label{eq:defn_Z}
			Z(\ell) \coloneqq \bigcup_k \{T_k-T_0=\ell \} \: \text{ and } \: 
			\widehat{Z}_{M-k}(\ell) \coloneqq Z(\ell) \cap \bigcap_{j=\ell +1}^{M-k}(Z(j))^\compl. 
		\end{align}
		Here, $Z(\ell)$ is the event that there is a $k$ such that $T_k-T_0=\ell$, i.e., a regeneration occurs after exactly $\ell$ steps, and $\widehat{Z}_{M-k}(\ell)$ is the event that in addition, then no further regeneration occurs until the $(M-k)$-th step.
		
		\begin{lemma}[Lemma~A.4 in \cite{SteibersPhD2017}]
			\label{lem:ADEregen}
			Let $d\ge 1$. Fix $M \in [2N/5, N]$ 
			and some initial position $z=(y,m)\in \widetilde{\cP}(N)$.
			Then the following holds:
			\begin{itemize}
				\item[i)] For every $\ell \le M-m$ and $x\in\bZ^d$
				\begin{equation}
					\label{eq:lemmaA4_1}
					\bP^z\big( X_{m+\ell} =x \,\vert\, \widehat{Z}_{M-m}(\ell) \big) \le C\ell^{-\frac{d}{2}}.
				\end{equation}
				\item[ii)] For every $\ell \le M-m$, $x\in\bZ^d$ and all unit vectors $e_j\in\bZ^d$
				\begin{equation}
					\label{eq:lemmaA4_2}
					\Abs{\bP^{(y,m)}\big( X_{m+\ell}=x\,\vert\, \widehat{Z}_{M-m}(\ell) \big) - \bP^{(y+e_j,m)}\big( X_{m+\ell}=x \,\vert\, \widehat{Z}_{M-m}(\ell) \big)} \le C\ell^{-\frac{d+1}{2}}.
				\end{equation}
				\item[iii)] For every $\ell\le M-m$ and $x\in\bZ^d$
				\begin{equation}
					\label{eq:lemmaA4_3}
					\Abs{\bP^{(y,m)}\big( X_{m+\ell}=x \,\vert\, \widehat{Z}_{M-m}(\ell) \big)-\bP^{(y,m+1)}\big( X_{m+\ell}=x \,\vert\, \widehat{Z}_{M-m-1}(\ell -1) \big)} \le C\ell^{-\frac{d+1}{2}}.
				\end{equation}
			\end{itemize}
		\end{lemma}
		
		We now provide a brief sketch of the proof.
		
		\begin{proof}
			For  \eqref{eq:lemmaA4_1} we will go into more detail since this is also the general approach for the other two estimates. First note that
			\begin{equation}
				\bP^z(Z(\ell)) \mathop{\longrightarrow}_{\ell \to \infty} \frac{1}{\mu_T}  \qquad
				\text{with } \;\; \mu_T \coloneqq \bE^z[T_2-T_1] \in (1,\infty)
			\end{equation}
			by the renewal theorem. (Recall that in our setting we in principle allow the random walk to start outside the percolation cluster, which possibly skews the law of the first increment $T_1-T_0$. However, by \cite[Lemma~B.1]{BethuelsenBirknerDepperschmidtSchluter}, 
			see also Lemma \ref{lem:quenched-rw-cluster} below,  
			the time until the path hits the cluster has exponential tails and can thus be neglected in the following arguments). For an intuitive idea why \eqref{eq:lemmaA4_1} holds observe that
			\[
			\bP^z(\{ X_{m+\ell}=x \} \cap Z(\ell)) = \sum_{k=1}^\ell \bP^z((X_{T_k},T_k)=(x,m+\ell))
			\]
			and by the CLT, there are $O(\sqrt{\ell})$ ``substantial'' terms in this sum (where $k = \ell/\mu_T \pm O(\sqrt{\ell})$), each of which is $O(\ell^{-(d+1)/2})$ by \eqref{eq:lemmaA2_1}.
			
			Here is the more detailed argument: Conditioned on $Z(\ell)$ the event $\{ X_{m+\ell}=x \}$ is independent of $\bigcap_{j=\ell +1}^{M-k}(Z(j))^\compl$ and thus we get for any $L\in\bN$
			\begin{align}
				\notag
				\bP^z\big( &X_{m+\ell} =x \,\vert\, \widehat{Z}_{M-m}(\ell) \big)=\bP^z(X_{m+\ell}=x \,\vert\, Z(\ell))\\
				\notag
				&= \frac{1}{\bP^z(Z(\ell))} \sum_{k=1}^L \bP^z((X_{T_k},T_k)=(x,m+\ell),T_{\frac{k}{2}}\ge \tfrac{\ell}{2}+m)\\
				\notag
				&\hspace{0.5cm}+ \frac{1}{\bP^z(Z(\ell))} \sum_{k=1}^L \bP^z((X_{T_k},T_k)=(x,m+\ell),T_k-T_{\frac{k}{2}}> \tfrac{\ell}{2})\\
				\notag
				&\hspace{0.5cm}+\frac{1}{\bP^z(Z(\ell))} \sum_{k=L+1}^\ell \bP^z((X_{T_k},T_k)=(x,m+\ell),T_{\frac{k}{2}}\le \tfrac{\ell}{2}+m)\\
				\label{eq:lemmaA4_proof_1}
				&\hspace{0.5cm}+ \frac{1}{\bP^z(Z(\ell))} \sum_{k=L+1}^\ell \bP^z((X_{T_k},T_k)=(x,m+\ell),T_k-T_{\frac{k}{2}}< \tfrac{\ell}{2}),
			\end{align}
			where for the second and fourth line we note that the events $\{T_{\frac{k}{2}}< \frac{\ell}{2}+m\}$ and $\{T_k-T_{\frac{k}{2}}>\frac{\ell}{2}\}$ coincide on $\{T_k=m+\ell\}$.
			
			Summing over all possible values of $(X_{T_{\frac{k}{2}}},T_{\frac{k}{2}})$, using \eqref{eq:lemmaA2_1} and translation invariance (and the regeneration property at time $T_{\frac{k}{2}}$) we obtain
			\begin{equation*}
				\bP^z\big((X_{T_k},T_k)=(x,m+\ell),T_{\frac{k}{2}}\le \tfrac{\ell}{2}+m\big) \le C k^{-\frac{d+1}{2}} \bP^z\big( T_{\frac{k}{2}}\le \tfrac{\ell}{2}+m\big).
			\end{equation*}
			Analogously we obtain similar bounds for the other terms, that is
			\begin{align*}
				\bP^z\big((X_{T_k},T_k)=(x,m+\ell),T_k-T_{\frac{k}{2}}<\tfrac{\ell}{2}\big) &\le Ck^{-\frac{d+1}{2}}\bP^z\big(T_{\frac{k}{2}}\le \tfrac{\ell}{2} +m\big),\\
				\bP^z\big((X_{T_k},T_k)=(x,m+\ell),T_{\frac{k}{2}}\ge \tfrac{\ell}{2}+m\big) &\le Ck^{-\frac{d+1}{2}}\bP^z\big(T_{\frac{k}{2}}\ge \tfrac{\ell}{2}+m\big),\\
				\bP^z\big((X_{T_k},T_k)=(x,m+\ell),T_k-T_{\frac{k}{2}}> \tfrac{\ell}{2}\big) &\le Ck^{-\frac{d+1}{2}}\bP^z\big(T_{\frac{k}{2}}\ge \tfrac{\ell}{2}+m\big).
			\end{align*}
			
			Setting 
			$L \coloneqq \ell/\mu_T$ we get, using the Markov-inequality and standard estimates on the $2d$-th moment of sums of centred i.i.d.\ random variables, e.g.\ \cite[Ch.~II.4, (PSF) on
			p.~152]{Perkins1999} or a Rosenthal inequality,
			\begin{align*}
				& \bP^z\big(T_{\frac{k}{2}} \ge \tfrac{\ell}{2}+m\big)
				\le \frac{\bE^z\big[ (T_{\tfrac{k}{2}} - m - \tfrac{k}{2} \mu_T)^{2d} \big]}{(\tfrac{\ell}{2}-\tfrac{k}{2} \mu_T)^{2d}}
				\le C\frac{k^d}{(L-k)^{2d}} \qquad \text{for } k< L
			\end{align*}
			and in a similar fashion
			\begin{equation*}
				\bP^z\big( T_{\frac{k}{2}} \le \tfrac{\ell}{2}+ m \big) \le C\frac{k^d}{(L-k)^{2d}} \qquad \text{for } k> L.
			\end{equation*}
			
			Note that  $\tfrac{k^d}{(L-k)^{2d}} >1$ only for $k\in[L-\sqrt{L},L+\sqrt{L}]$, and thus we have the upper bound
			\begin{align}
				\bP^z&(X_{m+\ell}=x \,\vert\, Z(\ell)) \notag \\
				& \le \sum_{k=1}^L k^{-\frac{d+1}{2}} \min\big\{ 1,C\tfrac{k^d}{(L-k)^{2d}}\big\} + \sum_{k=L+1}^\ell k^{-\frac{d+1}{2}} \min\big\{ 1,C\tfrac{k^d}{(k-L)^{2d}}\big\} \notag \\
				&\le \sum_{k=1}^{L-\sqrt{L}} Ck^{-\frac{d+1}{2}}\frac{k^d}{(L-k)^{2d}} + \sum_{k=L-\sqrt{L}}^{L+\sqrt{L}} C k^{-\frac{d+1}{2}} + \sum_{k=L+\sqrt{L}}^\infty Ck^{-\frac{d+1}{2}}\frac{k^d}{(L-k)^{2d}} \notag \\
				&\le C\int_0^{L-\sqrt{L}}\frac{t^{\frac{d-1}{2}}}{(L-t)^{2d}}\,dt + \sum_{k=L-\sqrt{L}}^{L+\sqrt{L}} C k^{-\frac{d+1}{2}} + C\int_{L+\sqrt{L}}^{\infty} \frac{t^{\frac{d-1}{2}}}{(L-t)^{2d}}\,dt \notag \\
				&\le C\ell^{-\frac{d}{2}}. \label{eq:lemmaA4_1.1}
			\end{align}
			For the integrals the upper bound $\ell^{-\frac{d}{2}}$ can be obtained using the substitution $\varphi(t)=L-t$. 
			
			For \eqref{eq:lemmaA4_2} we use the same approach as in \eqref{eq:lemmaA4_proof_1}. Then, with translation invariance and \eqref{eq:lemmaA2_3} we obtain
			\begin{align*}
				&\Abs{\bP^{(y,m)}\big( X_{m+\ell}=x\,\vert\, \widehat{Z}_{M-m}(\ell) \big) - \bP^{(y+e_j,m)}\big( X_{m+\ell}=x \,\vert\, \widehat{Z}_{M-m}(\ell) \big)}\\
				&\le C\Big( \sum_{k=1}^L k^{-\frac{d+2}{2}} \bP^z(T_{\frac{k}{2}}\ge \tfrac{\ell}{2}+m) + \sum_{k=L+1}^\infty k^{-\frac{d+2}{2}}\bP^z(T_{\frac{k}{2}}\le \tfrac{\ell}{2}+m) \Big)\\
				&\le C \sum_{k=1}^\infty k^{-\frac{d+2}{2}}\min\big\{ 1,\tfrac{k^d}{(L-k)^{2d}} \big\}\\
				&\le C \ell^{-\frac{d+1}{2}}.
			\end{align*}
			
			Similarly, for \eqref{eq:lemmaA4_3}, recall the definition of $Z(\ell)$ and $\widehat{Z}_{M-m}(\ell)$ from \eqref{eq:defn_Z}. Noting that $Z(\ell)$ is independent from the starting position and using the triangle inequality we get
			\begin{align*}
				&\Abs{\bP^{(y,m)}\big( X_{m+\ell}=x \,\vert\, \widehat{Z}_{M-m}(\ell) \big)-\bP^{(y,m+1)}\big( X_{m+\ell}=x \,\vert\, \widehat{Z}_{M-m-1}(\ell -1) \big)}\\
				&\le \Abs{\frac{1}{\bP^{(y,m+1)}(Z(\ell))} - \frac{1}{\bP^{(y,m+1)}(Z(\ell-1))}}\sum_{k=1}^\infty \bP^{(y,m+1)}((X_{T_k},T_k)=(x,m+\ell+1))\\
				&\hspace{1em}+\frac{1}{\bP^{(y,m+1)}(Z(\ell-1))} \sum_{k=1}^\infty \Big| \bP^{(y,m+1)}((X_{T_k},T_k)=(x,m+\ell+1)) \\[-2ex]
				& \hspace{14em} - \bP^{(y,m+1)}((X_{T_k},T_k)=(x,m+\ell)) \Big|.
			\end{align*}
			The last term can be handled in the same way as we did above for \eqref{eq:lemmaA4_2} using \eqref{eq:lemmaA2_2} instead of \eqref{eq:lemmaA2_3}. The first term is bounded from above by
			\begin{equation*}
				C\Abs{ \bP^{(y,m+1)}(Z(\ell-1))-\bP^{(y,m+1)}(Z(\ell)) } \ell^{-\frac{d}{2}},
			\end{equation*}
			where we used the estimate from \eqref{eq:lemmaA4_1}.
			
			It remains to observe that 
			\begin{equation}
				\label{eq:A1}
				\Abs{ \bP^{(y,m+1)}(Z(\ell-1))-\bP^{(y,m+1)}(Z(\ell)) } \le C\ell^{-\frac{1}{2}}.
			\end{equation}
			For example, classical results on the speed of convergence in Blackwell's renewal theorem
			show that \eqref{eq:A1} in fact holds if $\bE^z[(T_2-T_1)^{3/2}] < \infty$ (which is implied by the exponential tails of the regeneration times), see \cite[Proposition~1]{Lindvall1979}.
		\end{proof}
		
		\begin{remark}
			In the proof of Lemma~\ref{lem:ADEregen} we used (only) that the inter-regeneration times
			and spatial increments have finite $2d$-th moments, with the idea of potential re-usability for other models. In the situation considered here, we in fact have exponential tail bounds and thus even an exponential upper bound holds in \eqref{eq:A1}  by Theorem~3.1 from \cite{Baxendale2005}, or see \cite{Kendall1959}. Similarly, the arguments leading to \eqref{eq:lemmaA4_1.1} could be shortened when using exponential tails.
		\end{remark}

		\begin{proof}[Proof of Lemma~\ref{lem:annealed_derivative_estimates}]
			For completeness sake we will also include an annealed estimate which was not specifically mentioned in Lemma~\ref{lem:annealed_derivative_estimates}, namely
			\begin{equation}
				\label{eq:ann_der_est_proof_1}
				\bP^{(y,m)}(X_{n+m}=x) \le Cn^{-d/2}.
			\end{equation}
			To prove \eqref{eq:ann_der_est_proof_1} we use \eqref{eq:lemmaA4_1}, translation invariance and the exponential tails of the regeneration times
			\begin{align*}
				&\bP^{(y,m)}(X_{n+m}=x) \\
				&= \sum_{\ell \le n} \bP^{(y,m)}\big( \widehat{Z}_n(\ell) \big)\sum_{w\in\bZ^d} \bP^{(y,m)}\big( X_{m+\ell} =w \,\vert\, \widehat{Z}_n(\ell) \big) \bP^{(y,m)}\big( X_{n+m} =x \,\vert\, \widehat{Z}_n(\ell), X_{m+\ell}=w \big)\\
				& \le \sum_{\ell \le n} C \bP^{(y,m)}\big( \widehat{Z}_n(\ell) \big) \ell^{-\frac{d}{2}}
				\sum_{v \in \bZ^d} \bP^{(o,m)}\big( X_{n+m} = v \,\vert\, \widehat{Z}_n(\ell), X_{m+\ell}= o \big) \\
				&\le \sum_{\ell \le n} C e^{-c(n-\ell)} \ell^{-\frac{d}{2}}\\
				&\le \sum_{\ell \le \frac{n}{2}} Ce^{-c n/2} \ell^{-\frac{d}{2}} + \sum_{ \frac{n}{2} \le \ell \le n} Ce^{-c (n-\ell)} \Big(\frac{n}{2} \Big)^{-\frac{d}{2}}\\
				&\le C n^{-d/2}.
			\end{align*}
			Next we prove \eqref{eq:1} and \eqref{eq:3}. For \eqref{eq:4} and \eqref{eq:5} we note that they then follow by translation invariance. Using the ideas from above, specifically changing the summation for $\bP^{(y,m)}\big( X_{n+m} =x \,\vert\, \widehat{Z}_n(\ell), X_{m+\ell}=w \big)$ as was done in the last display, the triangle inequality and \eqref{eq:lemmaA4_2} we get the estimate for \eqref{eq:1}:
			\begin{align*}
				&\Abs{ \bP^{(y,m)}(X_{n+m}=x) - \bP^{(y+e_j,m)}(X_{n+m}=x) }\\
				&\begin{multlined}[t]
					\le \sum_{\ell \le n} \bP^{(y,m)}(\widehat{Z}_n(\ell))\\
					\hspace{1cm} \cdot \sum_{w\in\bZ^d} \Abs{\bP^{(y,m)}\big(X_{m+\ell}=w \,\vert\, \widehat{Z}_n(\ell)\big)- \bP^{(y,m)}\big(X_{m+\ell}=w-e_j \,\vert\, \widehat{Z}_n(\ell)\big)}\\
					\cdot\bP^{(y,m)}\big( X_{n+m}=x \,\vert\, \widehat{Z}_n(\ell), X_{m+\ell}=w \big)
				\end{multlined}\\
				&\le \sum_{\ell \le n} Ce^{-c(n-\ell)}\ell^{-\frac{d+1}{2}}\\
				&\le Cn^{-\frac{d+1}{2}}.
			\end{align*}
			For \eqref{eq:3} we can use the same ideas to show
			\begin{align*}
				&\Abs{\bP^{(y,m)}(X_{n+m}=x) - \bP^{(y,m+1)}(X_{n+m}=x) }\\
				&\le \sum_{\ell \le n} \bP^{(y,m)}\big(\widehat{Z}_n(\ell)\big)\sum_{w\in\bZ^d}\Abs{\bP^{(y,m)}\big(X_{m+\ell}=w \,\vert\, \widehat{Z}_n(\ell) \big) -\bP^{(y,m+1)}\big( X_{m+\ell}=w \,\vert\, \widehat{Z}_{n-1}(\ell -1) \big)}\\
				& \hspace{8cm} \cdot \bP^{(y,m)}\big(X_{n+m}=x \,\vert\, X_{m+\ell}=w, \widehat{Z}_n(\ell)\big) \\
				&\begin{multlined}[t]
					+\sum_{\ell \le n} \Abs{ \bP^{(y,m)}\big( \widehat{Z}_n(\ell) \big) - \bP^{(y,m+1)}\big( \widehat{Z}_{n-1}(\ell -1) \big) }\\
					\cdot \sum_{w\in\bZ^d} \bP^{(y,m+1)}\big( X_{m+\ell}=w \,\vert\, \widehat{Z}_{n-1}(\ell -1) \big)\bP^{(y,m)}\big( X_{n+m}=x \,\vert\, X_{m+\ell}=w,\widehat{Z}_n(\ell) \big).
				\end{multlined}
			\end{align*}
			Note that the first sum can be treated in the same way we proved the upper bound for \eqref{eq:1} using \eqref{eq:lemmaA4_3} instead of \eqref{eq:lemmaA4_2}. For the second sum note that, due to translation invariance
			\begin{align*}
				&\Abs{ \bP^{(y,m)}\big( \widehat{Z}_n(\ell) \big) - \bP^{(y,m+1)}\big( \widehat{Z}_{n-1}(\ell -1) \big) }\\
				&=\Abs{\bP^{(y,m)}\big( Z(\ell) \big) - \bP^{(y,m+1)}\big( Z(\ell-1) \big)} \bP^{(y,m)}\big( \bigcap_{j=\ell +1}^n Z(j)^\compl \,\vert\, Z(\ell) \big)\\
				&\le C\ell^{-\frac{1}{2}} e^{-c(n-\ell)},
			\end{align*}
			where we used the estimate \eqref{eq:A1} from the previous proof and exponential tail bounds of the regeneration times, see Lemma~2.5 in \cite{BirknerCernyDepperschmidtGantert2013}. Using this in the estimate above, \eqref{eq:lemmaA4_1} and translation invariance we obtain
			\begin{align*}
				&\Abs{\bP^{(y,m)}(X_{n+m}=x) - \bP^{(y,m+1)}(X_{n+m}=x) }\\
				&\le Cn^{-\frac{d+1}{2}}+\sum_{\ell \le n} C\ell^{-\frac{1}{2}} e^{-c(n-\ell)} \sum_{w\in\bZ^d} \bP^{(y,m+1)}\big( X_{m+\ell}=w \,\vert\, \widehat{Z}_{n-1}(\ell -1) \big)\\
				&\hspace{6cm}\cdot \bP^{(y,m)}\big( X_{n+m}=x \,\vert\, X_{m+\ell}=w,\widehat{Z}_n(\ell) \big)\\
				&\le Cn^{-\frac{d+1}{2}}+\sum_{\ell \le n} C\ell^{-\frac{1}{2}} e^{-c(n-\ell)} \sum_{w\in\bZ^d} \ell^{-\frac{d}{2}} \bP^{(y,m)}\big( X_{n+m}=x \,\vert\, X_{m+\ell}=w,\widehat{Z}_n(\ell) \big)\\
				&\le Cn^{-\frac{d+1}{2}}+\sum_{\ell \le n} C\ell^{-\frac{d+1}{2}} e^{-c(n-\ell)}\\
				&\le Cn^{-\frac{d+1}{2}}.
			\end{align*}
			As mentioned above the other two estimates can then be obtained by translation invariance which concludes the proof.
		\end{proof}
		
		\subsection{Proof of Lemma \ref{lem:additional_annealed_estimate}}
		The fact that Lemma~\ref{lem:additional_annealed_estimate} holds for all $d\ge 1$ is a direct consequence of the extension to all $d\ge 1$ in Lemma~\ref{lem:lem3.6SS} and Lemma~\ref{lem:annealed_derivative_estimates}. For completeness sake we briefly discuss the proof here. Consider the following set of boxes centred at the origin of $\bZ^d$:
		\begin{align*}
			\widetilde{\Pi}^{(\varepsilon)}_n\coloneqq \{
			B \in\Pi^{(\varepsilon)}_n \colon B \cap[-\sqrt{n}\log^3
			n,\sqrt{n}\log^3n]^d\neq \emptyset \}.
		\end{align*}
		With this notation we can write the sum on the left hand side of
		\eqref{eq:6} as
		\begin{align}
			\label{eq:proof_extra_annealed_estimations_1}
			& \sum_{ B\in \widetilde{\Pi}^{(\varepsilon)}_n}
			\sum_{x\in B} \max_{y\in B}\bigl[\bP^{(o,0)}(X_n=y)-\bP^{(o,0)}(X_n=x)\bigr]\\
			\label{eq:proof_extra_annealed_estimations_2}
			& \qquad
			+ \sum_{B \in \Pi^{(\varepsilon)}_n\setminus\widetilde{\Pi}^{(\varepsilon)}_n}
			\sum_{x\in B}
			\max_{y\in B}\bigl[\bP^{(o,0)}(X_n=y)-\bP^{(o,0)}(X_n=x)\bigr].
		\end{align}
		It suffices to prove suitable upper bounds for these two sums.
		By Lemma~\ref{lem:lem3.6SS} we have
		\begin{align}
			\label{eq:proof_extra_annealed_estimations_3}
			\sum_{B \in \Pi^{(\varepsilon)}_n\setminus\widetilde{\Pi}^{(\varepsilon)}_n}
			\bP^{(o,0)}(X_n\in B) \le Cn^{-c\log n}
		\end{align}
		for some positive constants $C$ and $c$. Thus, the double sum
		\eqref{eq:proof_extra_annealed_estimations_2} is bounded from above
		by
		\begin{multline*}
			\sum_{B \in \Pi^{(\varepsilon)}_n\setminus\widetilde{\Pi}^{(\varepsilon)}_n}
			\sum_{x\in B} \bigl[\bP^{(o,0)}(X_n\in B) - \bP^{(o,0)}(X_n=x)\bigr] \\
			= \sum_{B \in \Pi^{(\varepsilon)}_n\setminus\widetilde{\Pi}^{(\varepsilon)}_n}
			(\abs{B}-1)\bP^{(o,0)}(X_n\in B)
			\le Cn^{d\varepsilon}n^{-c\log n} \le \widetilde{C}n^{-\tilde{c}\log n}
		\end{multline*}
		for suitably chosen constants $\tilde c$ and $\widetilde C$. Using
		annealed derivative estimates from
		Lemma~\ref{lem:annealed_derivative_estimates} the double sum
		\eqref{eq:proof_extra_annealed_estimations_1} is bounded above by
		\begin{align*}
			\sum_{B \in \widetilde{\Pi}^{(\varepsilon)}_n}\sum_{x\in B} Cn^\varepsilon
			n^{-\frac{d+1}{2}} \le C(n^\varepsilon+\sqrt{n}\log^3 n)^d
			n^\varepsilon n^{-\frac{d+1}{2}} \le Cn^{3d\varepsilon}n^{-1/2}.
		\end{align*}
		Combination of the last two displays completes the proof.
		
		\subsection{Proof of Lemma \ref{lem:lem3.6SS}}
		\label{sec:pf_lem3.6SS}
		The estimate in \eqref{eq:9} is a direct consequence of the exponential tail bounds of the regeneration times. We follow the same ideas as in the proof of Lemma~\ref{lem:Prop 3.11 Ersatz}. To make this section more accessible we repeat these arguments.
		
		Recall the event $R^\simu_N$ from \eqref{eq:defn_RN} and note that $\bP(R_N^\simu) \ge 1-CN^{-c\log N}$, again due to the exponential tail bounds. By translation invariance, we may restrict the proof to $(y,m)=(o,0)$. On the event $R^\simu_N$ the time between two consecutive regenerations is at most $\log^2 N$ up to time $N$. Moreover $(X_{T_k}-X_{T_{k-1}}, T_k-T_{k-1})_{1\le k\le N}$ conditioned on $R_N^\simu$ is an i.i.d.\ sequence and the increments $X_{T_k}-X_{T_{k-1}}$ are symmetrically distributed. Therefore
		\begin{align*}
			\bP^{(o,0)}\big( \norm{X_n} \ge \sqrt{n}\log^3 (N) \big)&\le \bP^{(o,0)}\big( \norm{X_n} \ge \sqrt{n}\log^3(N) \,\vert\, R^\simu_N \big) + CN^{-c\log N}\\
			&\le \bP^{(o,0)}\big( \exists k\le n \colon \norm{X_{T_k}}\ge \frac{1}{2}\sqrt{n}\log^3(N)\,\vert\, R^\simu_N \big) + CN^{-c\log N}\\
			&\le \sum_{k=1}^n \bP^{(o,0)}\big( \norm{X_{T_k}}\ge \frac{1}{2}\sqrt{n}\log^3(N)\,\vert\, R^\simu_N \big) + CN^{-c\log N}\\
			&\le d\sum_{k=1}^n \exp\bigg( -C\frac{n\log^6(N)}{k\log^4(N)} \bigg) +CN^{-c\log N}\\
			&\le CN^{-c\log N},
		\end{align*}
		where we used the fact that on $R^\simu_N$ the random walk cannot reach the distance $\sqrt{n}\log^3(N)$ and come back to $\tfrac{1}{2}\sqrt{n}\log^3(N)$ in between two regeneration times and Azuma's inequality for each coordinate in the fourth line. This proves \eqref{eq:9}. 
		
		Lastly, note that by \eqref{eq:9} and using the Markov inequality
		\begin{align*}
			\bP&\bigg( \Big\{ \omega\in\Omega \colon P^{(o,0)}_\omega\big(\norm{X_n}\ge \sqrt{n}\log^3(N) \big) \ge \sqrt{CN^{-c\log N}}\Big\} \bigg)\\
			&\le \frac{\bE\Big[ P^{(o,0)}_\omega\big(\norm{X_n}\ge \sqrt{n}\log^3(N)\big) \Big]}{\sqrt{CN^{-c\log N}}}\\
			&\le CN^{-\frac{c}{2}\log N},
		\end{align*}
		i.e.\  \eqref{eq:10} holds too.
		
		\section{On Lemma~B.1 from \cite{BethuelsenBirknerDepperschmidtSchluter}
			and Lemma~2.11 from \cite{BirknerCernyDepperschmidtRWDRE2016}}
		
		In Lemma~B.1 of \cite{BethuelsenBirknerDepperschmidtSchluter}, the following bound was
		stated and (partly) proved.         
		\begin{lemma}
			Let \label{lem:quenched-rw-cluster} $d\geq 1$ and, for  $C',c' \in (0,\infty)$, denote by $A_n(C',c')$ the set
			\begin{equation}
				\{ \omega\in\Omega:
				P^{(0,0)}_\omega(\xi_i(X_i)=0,i=1,\dots,n)\leq C'\mathrm{e}^{-c'n}
				\}.\end{equation}
			For $p>p_c(d)$ there are constants $C,c \in (0,\infty)$ so that 
			\begin{align*}
				\bP(A^\compl_n) \leq C\mathrm{e}^{-cn} \quad \text{for all $n=1,2,\dots$},
			\end{align*}
			where $A_n=A_n(C',c')$ for some small enough $C'$ and $c'$.
		\end{lemma}

		It has been brought to our attention that the proof given in
		\cite{BethuelsenBirknerDepperschmidtSchluter}, which used
		\cite[Lemma~2.11]{BirknerCernyDepperschmidtRWDRE2016}, might
		only work for $p$ sufficiently close to $1$ because the application of
		\cite[Lemma~2.11]{BirknerCernyDepperschmidtRWDRE2016}
		implicitly assumes this.  In fact, the crucial inequality from
		\cite[Lemma~2.11]{BirknerCernyDepperschmidtRWDRE2016} holds
		for any $d\geq 1$ and any $p>p_c$. In our present notation,
		this reads as follows.
		\begin{lemma}[\protect{Extension  of \cite[Lemma~2.11]{BirknerCernyDepperschmidtRWDRE2016}}]
			For \label{lem:drysitesbd} any $p>p_c$ there exist
			$\widetilde{C}<\infty$ and $\widetilde{c}>0$ such that for
			any
			$V = \{ (x_i, t_i) : 1 \le i \le k \} \subset \bZ^d \times
			\bZ$ with $0\leq t_1 < t_2 < \cdots < t_k$, we have
			\begin{align}
				\label{eq:drysitesbd}
				\bP\bigl( \xi_t(x)=0 \; \text{for all}\; (x,t) \in V \bigr) \leq \widetilde{C} e^{-\widetilde{c} k}.
			\end{align}
			
		\end{lemma}
		\begin{proof} We will argue analogously to
			\cite[Lemma~2.11]{BirknerCernyDepperschmidtRWDRE2016} but
			use a more precise estimate. The general idea is that the
			event in \eqref{eq:drysitesbd} either enforces the existence
			of many ($\ge \varepsilon k$) disjoint small clusters with
			specified starting points or the existence of a few large
			but finite clusters whose combined height exceed $k$; both
			types of events are unlikely for $p>p_c$.  Note that in
			contrast to the set-up in
			\cite{BirknerCernyDepperschmidtRWDRE2016}, we consider here
			the oriented cluster in the forward time direction.
			
			We assume w.l.o.g.\ that $t_i - t_{i-1} \geq 2$ for all $i$, as otherwise we may replace $V$ by a suitable
			subset of size $\geq \lfloor k/2 \rfloor$. 
			Consider the extension of the sequence $t_1 < \cdots < t_k$ to an infinite sequence by letting $t_{k+\ell} \coloneqq
			t_k+2\ell$ for $\ell \in \bN$. Moreover, define
			\[
			D(i) \coloneqq \inf \big\{ j > i : (x_i,t_i) \not\to^\omega \bZ^d \times \{ t_j-1 \} \big\},
			\quad i=1,2,\dots,k.
			\]
			Thus, $D(i)-1$ is the \emph{height} of the oriented cluster (when measured along the sequence of the $t_j$'s)
			which starts at $(x_i,t_i)$ and so $\xi_{t_i}(x_i)=0$ is equivalent to $D(i)<\infty$.
			By decomposing the clusters starting at the space-time points in $V$ into non-overlapping
			time slices (see \cite{BirknerCernyDepperschmidtRWDRE2016}, proof of Lemma~2.11 for a more
			detailed description and also Figure~1 there), we see that 
			\[
			\big\{ D(1)<\infty,\dots,D(k)<\infty\big\} = \bigcup_{m=1}^k A_{k,m}
			\]
			with (we implicitly put $d_0 = 0$)
			\begin{align*}
				A_{k,m} \coloneqq \bigcup_{1<d_1<\cdots<d_{m-1} \le k < d_m}
				\big\{ D(d_{i-1}) = d_i \text{ for }i=1,\dots,m \big\} .
			\end{align*}
			Note that the event $\{ D(d_{i-1})=d_i \}$ is measurable w.r.t.\
			$\sigma( \omega(x,t) : x\in\bZ^d, t_{d_{i-1}} \le t < t_{d_i})$ and thus
			for fixed $d_1<\cdots<d_m$, the events $D(d_{i-1})=d_i$, $i=1,2,\dots,m$ are independent. 
			Particularly, by translation invariance, it holds that
			\begin{align*}
				&\bP(D(d_{i-1}) = d_i) \le \bP(d_i-d_{i-1} \le D(0) < \infty)
				\le C e^{-c (d_i-d_{i-1})},
			\end{align*}
			where $D(0)=\inf \big\{ j>0 : (o,0) \not\to^\omega \bZ^d \times \{ j \} \big\}$ and
			where, for the second inequality, we have used a well-known bound on the height of clusters in supercritical oriented percolation
			(see e.g.\ \cite[Lemma~A.1]{BirknerCernyDepperschmidtGantert2013} for details).   
			This gives
			\begin{align}
				\bP(A_{k,m})
				& \le \sum_{1<d_1<\cdots<d_{m-1} \le k < d_m} C^k  e^{-c d_m} \notag \\
				& = C^k e^{-c k} \sum_{1<d_1<\cdots<d_{m-1} \le k} \sum_{\ell=0}^\infty e^{-c \ell}
				\leq \frac{C^k e^{-c k}}{1-e^{-c}} \binom{k}{m-1}
				\label{eq:L2.11repl.bd1}
			\end{align}
			(and in fact only this type of estimate was used in
			\cite{BirknerCernyDepperschmidtRWDRE2016}).
			
			On the other hand, for fixed $d_1<\cdots<d_{m-1}$, the
			events $\{ D(d_{m-1})=d_m \}$ are disjoint for different
			choices of $d_m$ ($>d_{m-1}$), hence
			\begin{align*}
				& \bP(A_{k,m}) \\
				& \,
				= \bP\left( \bigcup_{d_m=k}^\infty \bigg( \bigg( \bigcup_{1<d_1<\cdots<d_{m-1} \le k} \big\{ D(d_{i-1}) = d_i \text{ for }i=1,\dots,m-1 \big\}
				\bigg)
				\cap \big\{ D(d_{m-1}) = d_m\big\} \bigg) \right) \\
				& \, \leq
				\bP\left( \bigcup_{1<d_1<\cdots<d_{m-1}} \big\{ D(d_{i-1}) = d_i \text{ for }i=1,\dots,m-1 \big\}
				\right) \bP(\xi_0(0)=0),
			\end{align*}
			and iterating this gives
			\begin{align}
				\label{eq:L2.11repl.bd2}
				\bP(A_{k,m}) \le \bP(\xi_0(0)=0)^m.
			\end{align}
			
			Now, combining \eqref{eq:L2.11repl.bd1} and \eqref{eq:L2.11repl.bd2}, we find that, for $\varepsilon \in (0,1)$ (which we suitable tune below)
			\begin{align}
				\bP(\xi_{t_i}(x_i)=0 \text{ for } 1 \le i \le k) 
				& \leq \sum_{m=1}^k \bP(A_{k,m}) \notag \\
				& \leq \frac{C^k e^{- c k}}{1-e^{-c}}
				\sum_{m=1}^{\lfloor \varepsilon k \rfloor} \binom{k}{m-1}
				+ \sum_{m=\lfloor \varepsilon k \rfloor+1}^\infty \bP(\xi_0(0)=0)^m.
				\label{eq:L2.11repl.bd3}
			\end{align}
			Moreover, using the upper bound in Cram\'er's theorem (see, e.g., \cite{DemboZeitouni},
			Thm.~2.2.3 and its proof), we have
			\begin{align*}
				C^{\varepsilon k} e^{-c k} \sum_{m=1}^{\lfloor \varepsilon k \rfloor} \binom{k}{m-1}
				& = \exp\big(-(c -\varepsilon \log(C))k\big)
				\, 2^k \mathrm{Bin}_{k,1/2}\big([0,\lfloor \varepsilon k \rfloor-1] \cap \bZ) \\[-1ex]
				& \leq \exp\Big( - k\big( c - \varepsilon\log(C)
				- I_{1/2}(\varepsilon) + \log(2) \big)\Big) \\
				& = \exp\Big(- k\big( c - \varepsilon\log(C)
				- \varepsilon\log(\varepsilon) - (1-\varepsilon)\log(1-\varepsilon)\big)\Big),
			\end{align*}
			where $I_{1/2}(x) = x\log(x) + (1-x)\log(1-x) -\log(2)$ is the rate function of
			the binomial family $\mathrm{Bin}_{k,p}$ with $p=1/2$.
			We choose $\varepsilon>0$ small enough so that the term on the right-hand side decays
			exponentially in $k$ and note that the second term on the right-hand side of \eqref{eq:L2.11repl.bd3}
			also decays exponentially in $k$ (for any choice of $\varepsilon>0$) to complete the proof.
		\end{proof}
		
		Now for any $p>p_c$, Lemma~\ref{lem:quenched-rw-cluster} can be proved
		exactly as in \cite[Appendix~B]{BethuelsenBirknerDepperschmidtSchluter}
		by replacing the application of \cite[Lemma~2.11]{BirknerCernyDepperschmidtRWDRE2016}
		there by that of Lemma~\ref{lem:drysitesbd} above. For ease of reference, we provide the
		short proof here.
		\begin{proof}[Proof of Lemma~\ref{lem:quenched-rw-cluster}]
			Note that by our definition of the quenched law, see
			equation~\eqref{eq:defn_quenched_law}, the quenched random
			walk performs a simple random walk until it hits the cluster
			$\mathcal{C}$. Thus, on the event that the random walk does not
			hit the cluster, we can replace the random walk with a simple
			random walk $(Y_n)_n$ that is independent of the
			environment. Therefore, we obtain that 
			\begin{align*}
				\bP^{(0,0)}\big(\xi_{0}(X_0)
				& =\dots=\xi_n(X_n)=0\big)\\
				& = \sum_{x_1,\dots,x_n} \bP^{(0,0)}\big((X_1,\dots,X_n)=(x_1,\dots,x_n),
				\xi_0(0)=\dots=\xi_n(x_n)=0 \big)\\
				& = \sum_{x_1,\dots,x_n} \bP^{(0,0)}\big((Y_1,\dots,Y_n)=(x_1,\dots,x_n),
				\xi_0(0)=\dots=\xi_n(x_n)=0 \big)\\
				& = \sum_{x_1,\dots,x_n} \bP^{(0,0)}\big((Y_1,\dots,Y_n)=(x_1,\dots,x_n)\big)
				\bP\big(\xi_0(0)=\dots=\xi_n(x_n)=0 \big)\\
				& \leq \widetilde{C} \mathrm{e}^{-\widetilde{c} n},
			\end{align*}
			where $\widetilde{C} $ and $\widetilde{c}$ are the constants from Lemma~\ref{lem:drysitesbd}, which depend only
			on $d$ and $p$\,($>p_c$). 
			Using this and the definition of the annealed law we therefore get that 
			\begin{align*}
				\widetilde{C} \mathrm{e}^{-\widetilde{c} n}
				& \geq \bP^{(0,0)}\big(\xi_{0}(X_0)=\dots=\xi_n(X_n)=0\big)\\
				& =\int_{A_n}
				P^{(0,0)}_\omega(\xi_i(X_i)=0,i=1,\dots,n) \,d
				\bP(\omega) \\
				& \qquad \qquad + \int_{A^\compl_n} P^{(0,0)}_\omega(\xi_i(X_i)=0,i=1,\dots,n) \, d\bP(\omega)\\
				& \geq \int_{A^\compl_n} P^{(0,0)}_\omega(\xi_i(X_i)=0,i=1,\dots,n) \, d\bP(\omega)\\
				& \geq \bP(A^\compl_n)C'\mathrm{e}^{-c'n}.
			\end{align*}
			Consequently,  we obtain that
			$\bP(A_n^\compl)\leq (\widetilde{C}/C') \mathrm{e}^{-cn}$ with
			$c=\widetilde{c} - c'>0$ by choosing $c'<\widetilde{c}$.
		\end{proof}
		
	\end{appendix} 
	
	\subsection*{Acknowledgement} \noindent We would like to thank two
	anonymous referees whose careful reading and insightful comments on a
	previous version helped to remove inaccuracies and to improve the
	presentation.


	
	

\end{document}